\documentclass[12pt]{article}

\RequirePackage{amsthm,amsmath,amsfonts,amssymb}
\RequirePackage{natbib}

\RequirePackage[colorlinks,citecolor=blue,urlcolor=blue,backref=page,backref=page]{hyperref}

\usepackage{verbatim}
\RequirePackage{mathtools}
\usepackage{bm}

\usepackage{comment}
\usepackage[inline]{enumitem}

\theoremstyle{plain}

\usepackage{graphicx,booktabs,float}
\usepackage{tikz}
\usetikzlibrary{arrows.meta, positioning, calc, fit, shapes.geometric, backgrounds}
\usepackage{algorithm}
\usepackage{algpseudocode}

 \newcommand{\bA}{ {\mathbf{A}} }

 \newcommand{\bG}{ {\mathbf{G}} }
 
 \newcommand{\bH}{ {\mathbf{H}} }
 
 \newcommand{\bI}{ {\mathbf{I}} }

 \newcommand{\bS}{ {\mathbf{S}} }

 \newcommand{\bX}{ {\mathbf{X}} }
 \newcommand{\bx}{ {\mathbf{x}} }
 \newcommand{\bY}{ {\mathbf{Y}} }

 \newcommand{\bbeta}{ {\mathbf{\beta}} }

 \newcommand{\btheta}{ {\mathbf{\theta}} }

 \newcommand{\balpha}{ {\mathbf{\alpha}} }

\newcommand{\calA}{\mathcal{A}}

\usepackage[]{graphicx}
\chardef\bslash=`\\ % p. 424, TeXbook

\let\opt\pkg
\def\E{\textrm{E}}
\def\Var{\textrm{Var}}
\def\ps{\pi}
\def\balpha{\mathbf{\alpha}}
\def\balphahat{\widehat \mathbf{\alpha}}
\def\bbeta{\mathbf{\beta}}
\def\bbetahat{\widehat \mathbf{\beta}}

\makeatletter
\newif\ifctagsplit@
\newif\iftagsleft@
\def\gobble@tag{\@ifstar\@gobble\@gobble}
\def\invalid@tag#1{\@amsmath@err{#1}{\the\tag@help}\gobble@tag}
\def\dft@tag{\invalid@tag{\string\tag\space not allowed here}}
\def\default@tag{\let\tag\dft@tag}
\default@tag
\def\maketag@@{\@ifstar\maketag@@@\tagform@}
\def\maketag@@@#1{\hbox{\m@th\normalfont#1}}
\def\tagform@#1{\maketag@@@{(\ignorespaces#1\unskip)}}
\iftagsleft@
  \def\@eqnnum{\hbox to1sp{}\rlap{\normalfont\normalcolor
    \hskip -\displaywidth\tagform@\theequation}}
\else
  \def\@eqnnum{{\normalfont\normalcolor \tagform@\theequation}}
\fi
\def\thetag{\leavevmode\tagform@}
\def\make@df@tag{\@ifstar\make@df@tag@@\make@df@tag@@@}
\def\make@df@tag@@#1{%
  \gdef\df@tag{\maketag@@@{#1}\def\@currentlabel{#1}}}
\def\make@df@tag@@@#1{\gdef\df@tag{\tagform@{#1}%
  \toks@\@xp{\p@equation{#1}}\edef\@currentlabel{\the\toks@}}}
\let\ltx@label\label
\def\label@in@display{%
    \ifx\df@label\@empty\else
        \@amsmath@err{Multiple \string\label's:
            label '\df@label' will be lost}\@eha
    \fi
    \gdef\df@label
}
\def\eqref#1{\textup{\tagform@{\ref{#1}}}}
\makeatother

\def\E{\textrm{E}}
\def\Var{\textrm{Var}}
\def\ps{\pi}
\def\balpha{\mathbf{\alpha}}
\def\balphahat{\widehat \mathbf{\alpha}}
\def\bbeta{\mathbf{\beta}}
\def\bbetahat{\widehat \mathbf{\beta}}

\makeatletter
\newif\ifctagsplit@
\newif\iftagsleft@
\def\gobble@tag{\@ifstar\@gobble\@gobble}
\def\invalid@tag#1{\@amsmath@err{#1}{\the\tag@help}\gobble@tag}
\def\dft@tag{\invalid@tag{\string\tag\space not allowed here}}
\def\default@tag{\let\tag\dft@tag}
\default@tag
\def\maketag@@{\@ifstar\maketag@@@\tagform@}
\def\maketag@@@#1{\hbox{\m@th\normalfont#1}}
\def\tagform@#1{\maketag@@@{(\ignorespaces#1\unskip)}}
\iftagsleft@
  \def\@eqnnum{\hbox to1sp{}\rlap{\normalfont\normalcolor
    \hskip -\displaywidth\tagform@\theequation}}
\else
  \def\@eqnnum{{\normalfont\normalcolor \tagform@\theequation}}
\fi
\def\thetag{\leavevmode\tagform@}
\def\make@df@tag{\@ifstar\make@df@tag@@\make@df@tag@@@}
\def\make@df@tag@@#1{%
  \gdef\df@tag{\maketag@@@{#1}\def\@currentlabel{#1}}}
\def\make@df@tag@@@#1{\gdef\df@tag{\tagform@{#1}%
  \toks@\@xp{\p@equation{#1}}\edef\@currentlabel{\the\toks@}}}
\let\ltx@label\label
\def\label@in@display{%
    \ifx\df@label\@empty\else
        \@amsmath@err{Multiple \string\label's:
            label '\df@label' will be lost}\@eha
    \fi
    \gdef\df@label
}
\def\eqref#1{\textup{\tagform@{\ref{#1}}}}

\let\ilimits@\displaylimits

\let\DOTSI\relax

\def\int{\DOTSI\intop\ilimits@}
\def\oint{\DOTSI\ointop\ilimits@}
\def\intkern@{\mkern-6mu\mathchoice{\mkern-3mu}{}{}{}}
\def\intdots@{\mathchoice{\@cdots}%
 {{\cdotp}\mkern1.5mu{\cdotp}\mkern1.5mu{\cdotp}}%
 {{\cdotp}\mkern1mu{\cdotp}\mkern1mu{\cdotp}}%
 {{\cdotp}\mkern1mu{\cdotp}\mkern1mu{\cdotp}}}
\def\iint{\DOTSI\protect\ints@\tw@}
\def\iiint{\DOTSI\protect\ints@\thr@@}
\def\iiiint{\DOTSI\protect\ints@{4}}
\def\idotsint{\DOTSI\protect\ints@\z@}
\def\ints@#1{%
  \mkern-7mu\mathchoice{\mkern-2mu}{}{}{}%
  \mathop{\mkern7mu\mathchoice{\mkern2mu}{}{}{}%
    \intop\ifnum#1=\z@\intdots@
    \else\intkern@\fi
    \ifnum#1>\tw@\intop\intkern@\fi
    \ifnum#1>\thr@@\intop\intkern@\fi
    \intop
  }\ilimits@
}
\makeatother

\def\Cov{\textrm{Cov}}

\def\bZDS{{\bm{Z}}}
\def\bA{{\boldsymbol{A}}}

\def\bH{\bm{H}}
\def\bI{\mathbf{I}}

\def\zerovec{\boldsymbol{0}}
\renewcommand{\dbinom}[2]{\left(\begin{array}{@{}c@{}}#1\\[3pt]#2\end{array}\right)}

\makeatletter
\newcommand{\raisemath}[1]{\mathpalette{\raisem@th{#1}}}
\newcommand{\raisem@th}[3]{\raisebox{#1}{$#2#3$}}
\makeatother

\def\opt{{\mbox{opt}}}

\def\calA{\mathcal{A}}

\newcommand{\upto}[2]{\overline{#2}_{#1}}
\newcommand{\after}[2]{\widetilde{#2}_{#1}}
\newcommand{\uptoopt}[2]{{\overline{#2}}_{#1}{\hspace*{-0.6em} \raisebox{0.925ex}{\scriptsize{$\textrm{opt}$}}} \; \;}
\newcommand{\afteropt}[2]{{\widetilde{#2}}_{#1}{\hspace*{-1.3em} \raisebox{0.925ex}{\scriptsize{$\textrm{opt}$}}} \; \;}

\def\Normal{\textrm{Normal}}

\newcommand\expit[1]{\textrm{expit}\left\{ #1 \right\}}

\newcommand\aopt[1]{a_{#1}{\hspace*{-0.35em} \raisebox{0.925ex}{\scriptsize{$\textrm{opt}$}}}}
\newcommand\mopt[1]{m_{#1}{\hspace*{-0.45em} \raisebox{0.925ex}{\scriptsize{$\textrm{opt}$}}}}
\newcommand\gammaopt[1]{\gamma_{#1}{\hspace*{-0.20em} \raisebox{0.925ex}{\scriptsize{$\textrm{opt}$}}}}

\newcommand\Yhatopt[2]{\widehat Y_{#1}{\hspace*{-0.20em} \raisebox{0.925ex}{\scriptsize{$(#2,\textrm{opt})$}}}}

\makeatletter
\renewcommand\subsubsection{\@startsection{subsubsection}{3}{\z@}%
                                     {-3.25ex\@plus -1ex \@minus -.2ex}%
                                     {-1.5ex \@plus -.2ex}%
                                     {\normalfont\bfseries}}
\makeatother

\newcounter{myexample}
\newenvironment{myexample}[1][]{\refstepcounter{myexample}\par\medskip
   \noindent \textbf{Example~\themyexample: #1} \rmfamily}{\medskip}

\usepackage{listings}
\def\bA{\mathbf{A}}
\def\bH{\mathbf{H}}

\def\bbetahat{\widehat{\bm{\beta}}}
\def\bbeta{\bm{\beta}}
\def\balphahat{\widehat{\bm{\alpha}}}
\def\balpha{{\bm{\alpha}}}

\def\btheta{{\bm{\theta}}}

\begin{document}

\title{G-estimation for Dynamic Treatment Regimes\thanks{This paper is an extended version of `G-estimation for dynamic treatment regimes in the longitudinal setting', paper 7 in \textit{Adaptive Treatment Strategies in Practice: Planning Trials and Analyzing Data for Personalized Medicine}, MR Kosorok and EEM Moodie, Society for Industrial and Applied Mathematics, 2015.}
\author{David A. Stephens\\[6pt] \texttt{david.stephens@mcgill.ca} \\[25pt]
	Department of Mathematics and Statistics, McGill University, Canada\\
	}
}
\date{}

\maketitle

\index{G-estimation|(}
\begin{abstract}
The method of G-estimation in the single and multiple interval setting is discussed, and applications to the design of optimal dynamic treatment regimes are given.   We consider cases of binary treatments, treatments lying in a finite set, or where treatment can take arbitrary continuous values.  Crucial in the construction is the idea that there is a baseline `zero' level of treatment against which responses to other treatments are compared.  For continuous treatments, the numerical value zero provides this baseline.  For discrete-valued treatments, we choose a baseline from amongst the list of treatment levels.  G-estimation for dynamic regimes proceeds by proposing expected (counterfactual) response models relative to the response at baseline treatment.
    
\end{abstract}
\section{G-estimation and semiparametric regression}

\subsection{Connection to E-estimation}\index{E-estimation|(}

In the seminal paper \cite{Robins1992}, Robins, Mark and Newey consider the following semiparametric regression specification for continuous response $Y$ in terms of exposure of interest $A$ in the presence of predictors $\bX$:
\begin{equation}
Y_i = A_i \theta + g(\bX_i) + \epsilon_i \label{eq:PLM}
\end{equation}
with $\epsilon_i | \bX_i, A_i, i=1,\ldots,n$ conditionally independent residual terms with zero mean and constant variance $\sigma^2$.  In this expression, the function $g(.)$ is a nonparametrically specified regression function that depends on $\bX$ but not on $A$. In \cite{Robins1992}, a finite linear basis function representation of $g(.)$ is used, specifically where
\[
g(\bx) = \beta_0 + \sum_{k=1}^K \bbeta_k b_k(\bx)
\]
yielding the linear model $Y_i = A_i \theta + \bX_i \bbeta + \epsilon_i$, $i=1\ldots,n$, or in vector form
\begin{equation*}
\bY = \bA \theta + \bX \bbeta + \epsilon.
\label{eq:BLM}
\end{equation*}
The ordinary least squares (OLS) estimator for $(\theta,\beta^\top)^\top$ is the solution of the equation
\[
\dbinom{\bA^\top}{\bX^\top} (\bY - \bA \theta - \bX \bbeta) = \zerovec.
\]
Using a two-step approach, we see that the OLS estimator of $\theta$ can be obtained by first writing
\[
\bY - \bA \theta = \bX \bbeta + \epsilon
\]
to yield the implicit estimator
\[
\bbetahat(\theta) = (\bX^\top \bX)^{-1} \bX^\top (\bY - \bA \theta).
\]
By substituting back into the previous equation, we see that
\[
\bA^\top (\bY - \bA \theta - \bH (\bY - \bA \theta)) = \bA^\top (\bI_n - \bH)(\bY - \bA \theta) = \zerovec
\]
where $\bH = \bX (\bX^\top \bX)^{-1}\bX^\top$, and  where $\bI_n$ is the $n \times n$ identity matrix.  In the usual way, this yields the solution
\[
\widehat \theta = (\bA^\top (\bI_n - \bH) \bA  )^{-1}  (\bA^\top (\bI_n - \bH)  ) \bY
\]
as $(\bI_n - \bH)$ is idempotent and symmetric.  Suppose now we consider the linear model that treats $\bX$ as the design matrix and $\bA$ as the outcome, that is, the model that seeks to explain variation in treatment using predictors contained in $\bX$.  This model can be written $\bA = \bX \balpha + \epsilon^\ast$, and by the usual OLS theory fitted values in this model are given by $\widehat \bA = \bH \bA$.  Hence $(\bI_n - \bH) \bA = \bA - \widehat\bA$, and thus
\begin{equation*}
\widehat \theta =  (\bA^\top (\bA - \widehat \bA)  )^{-1}  ( \bA - \widehat \bA  )^\top \bY
= \frac{\sum\limits_{i=1}^n Y_i (A_i - \widehat A_i)}{\sum\limits_{i=1}^n A_i (A_i - \widehat A_i)}.
 \label{eq:RMN-I}
\end{equation*}
From this form, we observe a useful facet of the model: that $\widehat \theta$ is a solution to the estimating equation
\begin{equation*}
(\bY - \bA \theta)^\top (\bA - \widehat \bA) = \sum_{i=1}^n (Y_i - A_i \theta) (A_i - \widehat A_i) = 0.
\end{equation*}
This estimating equation is unbiased, having expectation zero, provided that the conditional expectation model for $A$ as a function of $\bX$ is correctly specified, and that $Y$ and $A$ are conditionally independent given $\bX$.    This may motivate us to focus only on the treatment model (modelling of $A$ given $\bX$), rather than the outcome model (modelling of $Y$ given $A$ and $\bX$), as the former may be more straightforward to specify from expert knowledge. However, from the original least squares problem, we also have that $\widehat \theta$ is a solution to
\[
\sum_{i=1}^n (Y_i - A_i \theta - \bX_i \bbetahat) A_i  = 0
\]
and thus conclude also that we must have
\begin{equation*}
\sum_{i=1}^n (Y_i - A_i \widehat \theta - \bX_i \bbetahat) (A_i - \widehat A_i)  = 0 \label{eq:RMN-EE}
\end{equation*}
so that
\begin{equation}
\widehat \theta = \frac{\sum\limits_{i=1}^n (Y_i - \bX_i \bbetahat )(A_i - \widehat A_i)}{\sum\limits_{i=1}^n A_i (A_i - \widehat A_i)} \label{eq:RMN-II}
\end{equation}
as by construction  $\bX^\top  (\bA - \widehat \bA) = 0$, so
\[
(\bY - \bA \widehat \theta - \bX \bbetahat)^\top (\bA - \widehat \bA) = (\bY - \bA \widehat \theta )^\top (\bA - \widehat \bA) - \bbetahat^\top \bX^\top  (\bA - \widehat \bA) = 0.
\]
Robins et al. \cite{Robins1992} proceed to elaborate on the properties of general versions of the estimators \eqref{eq:RMN-I} and \eqref{eq:RMN-II} and demonstrate that they exhibit desirable asymptotic properties.  In particular they demonstrate that estimators of this form are consistent for the regression parameters and achieve the semiparametric variance bound provided the model for the dependence of $A$ on $X$ is correctly specified.  Robins et al. \cite{Robins1992} were not primarily concerned with causal inference\index{Causal inference}, but it transpires that the estimator in equation \eqref{eq:OLS1-feas} is consistent for $\theta$ provided the proposed conditional expectation model $\ps(\bX_i) = \E[A_i|\bX_i]$ is correctly specified: the estimator is singly robust.  A doubly robust\index{Doubly robust} extension is also given by \cite{Robins1992}; the estimator takes the form
\begin{equation}
\widehat \theta_{DR} = \frac{\sum\limits_{i=1}^n (Y_i-m_0(X_i; \bbetahat)) (A_i - \ps(\bX_i;\balphahat))}{\sum\limits_{i=1}^n A_i (A_i - \ps(\bX_i;\balphahat))}  \label{eq:OLS1-DR},
\end{equation}
where the function $m_0(X;\bbeta)$ is defined via
\[
m_0(X;\bbeta) = \E[Y|A,X] - \theta A = m(A,X) - \theta A
\]
say; the model $m(A,X)$ represents the conditional expectation of $Y$ given $(A,X)$. Parameter $\beta$ is estimated in a separate regression stage; this form matches that in \eqref{eq:RMN-II}.

For this particular model, where the dependence on the treatment $A$ is in this simple form of a linear term in $A$, the estimator in \eqref{eq:RMN-II} and \eqref{eq:OLS1-DR} is \emph{doubly robust}:\index{Doubly robust} it is consistent if at least one of $m_0(.;.)$ or $\ps(.;.)$ is correctly specified.   To be more specific, the corresponding doubly robust estimating equation provides a consistent estimator of the expectation of $Y$ as $A$ varies under the assumption that the entire effect of treatment $A$ is captured by the linear term in $A$; there are no higher order terms in $A$ (not possible in the binary treatment case) and no interaction terms.  The structural assumption that $\E[Y|A,X] - \theta A$ does not depend on $A$, that is,
\[
\E[Y|A,X] = \theta A + m_0(X;\bbeta) \qquad \textrm{a.s.}
\]
is critical to the doubly robust procedure.

\index{E-estimation|)}

\subsection{Extensions}

In the linear model case, the estimating equation implied by \eqref{eq:OLS1-DR} may be written
\begin{equation}
\sum_{i=1}^n (Y_i - m(A_i,X_i;\theta,\beta)) (A_i - \pi(X_i)) = 0
\label{eq:DREE}
\end{equation}
where $m(A_i,X_i;\theta,\beta) = m_0(X_i;\beta) + \gamma(A_i,X_i;\theta)$ is the conditional expectation of $Y_i$ given $(A_i,X_i)$, and $\ps(\bX_i) = \E[A_i|\bX_i]$ is the conditional expectation of $A_i$ given $X_i$.  This estimating equation is unbiased if at least one of the models $m(.;.)$ or $\pi(.)$ is correctly specified.

If $m_0(.;.)$  is identically zero, the estimator can be written
\begin{equation*}
\widehat \theta_{SR} = \frac{\sum\limits_{i=1}^n Y_i (A_i - \ps(\bX_i))}{\sum\limits_{i=1}^n A_i (A_i - \ps(\bX_i))}   \label{eq:OLS1-infeas}
\end{equation*}
which is the principal form given by \cite{Robins1992}, who demonstrate it to be the most efficient means of estimating $\theta$ in the model \eqref{eq:PLM} if function $g(.)$ is left unrestricted.  This infeasible `E-estimator' can be replaced by a feasible version by using an estimated version of $\ps(.)$, that is, we substitute a parametric model $\ps(.;\balpha)$ for $\pi(.)$, estimate parameters $\balpha$ in the first step, and construct the estimator
\begin{equation}
\widehat \theta_{SR} = \frac{\sum\limits_{i=1}^n Y_i (A_i - \ps(\bX_i;\balphahat))}{\sum\limits_{i=1}^n A_i (A_i - \ps(\bX_i;\balphahat))} .  \label{eq:OLS1-feas}
\end{equation}
Beyond the linear model formulation, in the setting where
\begin{equation}
\E[Y_i|A_i,X_i] = m_0(X_i;\beta) + \gamma(A_i,X_i;\btheta) \label{eq:GeneralRegressionForm}
\end{equation}
with $\gamma(A,X;\btheta)$ depending in a more general fashion on $A$, the OLS estimating equation becomes
\[
\sum_{i=1}^n (Y_i - m(A_i,\bX_i;\btheta,\bbeta)) \dbinom{T(\bX_i;\bbeta)}{S(A_i,\bX_i;\btheta)} = \zerovec
\]
where
\[
T(x;\bbeta) = \frac{\partial m_0(x;\bbeta)}{\partial \bbeta} \qquad \textrm{and} \qquad S(a,x;\btheta) = \frac{\partial \gamma(a,x;\btheta)}{\partial \btheta}
\]
are the column vectors of derivatives of the model components taken with respect to the parameters.  For example, if
\[
\gamma(a,x;\btheta) = a \theta_1 + a^2 \theta_2
\qquad \textrm{then} \qquad
\frac{\partial \gamma(a,x;\btheta)}{\partial \btheta} = \dbinom{a}{a^2};
\]
if
\[
\gamma(a,x;\btheta) = a \theta_1 + a x \theta_2 + a^2 \theta_3 + a^2 x \theta_4
\qquad \textrm{then} \qquad
\frac{\partial \gamma(a,x;\btheta)}{\partial \btheta} = \left(
\begin{array}{c}
  a \\
  a x \\
  a^2 \\
  a^2 x
\end{array}
\right)
\]
and so on. Following the logic above, and OLS principles, we construct the adjusted estimating equation
\begin{equation}
\sum_{i=1}^n (Y_i - m(A_i,\bX_i;\btheta,\beta)) \dbinom{T(\bX_i;\bbeta)}{S(A_i,\bX_i;\btheta)-\pi_S(\bX_i;\btheta)} = \zerovec
\label{eq:EE-S}
\end{equation}
where $\pi_S(X_i;\btheta) = \E_{A_i|X_i}[S(A_i,\bX_i;\btheta)|X_i]$.  This estimating equation is unbiased if at least one of the models encapsulated in $m_0$ or in $\pi_S$ is correctly specified and double robustness \index{Doubly robust} follows. From \eqref{eq:BLM}, the extension to multivariate treatments (applied in the single interval setting) is straightforward.  We merely allow the function $\gamma(A,X;\btheta)$ to depend on a vector-valued $\bA$, and proceed using identical estimating equations.

Note that a solution to estimating equations of the form of equation \eqref{eq:EE-S} can be obtained in a straightforward fashion using standard techniques whenever the outcome mean model $m(.,.;.)$ is linear in the parameters; in vector form the estimating equation is
\[
\bZDS^\top (\bY - \bX \varphi) =  \mathbf{0}
\]
where $\bX$ is the design matrix corresponding to the mean function, $\varphi = (\theta,\beta^\top)^\top$, and $\bZDS^\top$ is the matrix formed by binding the vectors
\[
\dbinom{T(\bX_i)}{S(A_i,\bX_i)-\pi_S(X_i)} \qquad i=1,\ldots,n
\]
columnwise.   We then have a closed form for the estimators as
\[
\widehat \varphi = (\bZDS^\top \bX)^{-1} \bZDS^\top \bY,
\]
akin to ordinary least squares.  More generally, if \eqref{eq:EE-S} is derived from a generalized linear model, iteratively reweighted least squares procedures can be used.  Finally, general iterative numerical procedures (e.g. Newton-like methods) can be used.

\subsection{G-estimation via one-step ordinary least squares}

In the simplest version of the model which assumes linear dependence on the parameters, we may rewrite \eqref{eq:GeneralRegressionForm} as 
\begin{align*}
\E[Y_i|A_i,\bX_i] & =  m_0(\bX_i;\bbeta) + \gamma(A_i,\bX_i;\btheta)  = m_0(\bX_i;\bbeta) + (\bG_i \odot \bS_i)^\top \btheta
\end{align*}
where
\[
(\bG_i \odot \bS_i) = \begin{bmatrix}
g_1(\bX_i) \\ g_2 (\bX_i) \\ \vdots \\ g_d(\bX_i) 
\end{bmatrix} \odot \begin{bmatrix}
s_1(A_i) \\ s_2 (A_i) \\ \vdots \\ s_d(A_i) 
\end{bmatrix}
\]
is an elementwise product. Then in the above notation
\[
\pi_S(\bX_i) = \E_{A_i|X_i}[(\bG_i \odot \bS_i)|\bX_i] = \bG_i \odot \E_{A_i|X_i}[\bS_i|\bX_i] = \bG_i \odot \pi_S(\bX_i)
\]
and we may rewrite the original form as
\[
\E[Y_i|A_i,\bX_i] =  m_0(\bX_i;\bbeta)+(\bG_i \odot \pi_S(\bX_i))^\top \btheta + (\bG_i \odot (\bS_i- \pi_S(\bX_i))^\top \btheta
\]
or equivalently
\begin{equation}\label{eq:AltGeneralRegressionForm}
\E[Y_i|A_i,\bX_i] =  \widetilde m_0(\bX_i;\widetilde \bbeta) + (\bG_i \odot \widetilde \bS_i)^\top \btheta
\end{equation}
where
\begin{align*}
\widetilde m_0(\bX_i;\widetilde \bbeta) & = m_0(\bX_i;\widetilde \bbeta)+(\bG_i \odot \pi_S(\bX_i))^\top \btheta \\[6pt]
\widetilde \bS_i & = \bS_i- \pi_S(\bX_i).
\end{align*}
Then it is clear that \eqref{eq:AltGeneralRegressionForm} can be used to fit the model directly using OLS, and that the resulting estimating equations match \eqref{eq:EE-S} without further adjustment.  For example, in the simplest case, if
\[
\E[Y_i|A_i,X_i] = \beta_0 + A_i \theta
\]
then we may rewrite this as
\[
\E[Y_i|A_i,X_i] = \widetilde m_0(\bX_i;\widetilde \bbeta) + (A_i-\pi(\bX_i)) \theta = 
\widetilde m_0(\bX_i;\widetilde \bbeta) + \widetilde S_i \theta
\]
where $\widetilde m_0(\bX_i;\widetilde \bbeta) = \beta_0 + \pi(\bX_i) \theta$.  Note that taking expectations with respect to $A_i$ given $\bX_i$ reveals that
\[
\E[Y_i|\bX_i] = \widetilde m_0(\bX_i;\widetilde \bbeta).
\]

\subsection{Connection with instrumental variables}
\index{Instrumental variables|(}
The estimator from \eqref{eq:OLS1-feas} can be viewed as a two-stage instrumental variable estimator.  The instrument is the derived residual $Z_i = A_i - \ps(\bX_i;\balpha)$; recall that for $Z_i$ to be an instrument for $A_i$, it
    \begin{enumerate}

    \item[(i)] is correlated with $A_i$, and

    \item[(ii)] is uncorrelated with $Y_i$, given $A_i$.

    \end{enumerate}
    For (i), we have that
\begin{eqnarray*}
\Cov[A_i,Z_i] = \Cov[A_i,A_i - \ps(\bX_i;\balpha)] &=& \E[A_i (A_i - \ps(\bX_i;\balpha)] = \E[A_i^2 - A_i \ps(\bX_i;\balpha)]\\[6pt]
& = & \E_{\bX_i} \left[ \E_{A_i|\bX_i}[ A_i^2 - A_i \E_{A_i|\bX_i}[ A_i | \bX_i]| \bX_i] \right]\\[6pt]
& = & \E_{\bX_i} \left[ \E_{A_i|\bX_i}[ A_i^2 | \bX_i] - \{\E_{A_i|\bX_i}[ A_i | \bX_i]\}^2\right]\\[6pt]
& = & \E_{\bX_i} \left[ \Var_{A_i|\bX_i}[ A_i| \bX_i]\right] > 0
\end{eqnarray*}
unless $\Var_{A_i|\bX_i}[ A_i| \bX_i] = 0$ with probability 1, in which case $\bX_i$ predicts $A_i$ precisely; beyond  that case, the covariance is positive.  For (ii) we have that
\begin{eqnarray*}
\E_{Y_i|A_i,Z_i}[Y_i|A_i,Z_i] & = & \E_{X_i|A_i,Z_i}\left[\E_{Y_i|A_i,Z_i,X_i}[Y_i|A_i,Z_i,X_i] \right]\\[6pt]
& \equiv & \E_{X_i|A_i}\left[\E_{Y_i|A_i,X_i}[Y_i|A_i,X_i]\right] \\[6pt]
& = & \E_{Y_i|A_i}[Y_i|A_i]
\end{eqnarray*}
as $Z_i$ is defined jointly via $(A_i,X_i)$.  Thus $Z_i$ is an instrument for $A_i$ in the relationship with $Y_i$.
\index{Instrumental variables|)}

\subsection{The general semiparametric formulation}

The OLS formulation can be generalized to the semiparametric setting, where the model space is the Hilbert space $\mathcal{H}_{A,X}$ of square integrable functions of $(A,X)$ equipped with an inner product defined by the covariance operator (see \cite{Tsiatis_Book}).

Informally, we assume that $m(A,X) \in \mathcal{H}_{A,X}$, and define the estimator by considering projections of data $Y$ onto $\mathcal{H}_{A,X}$.  Denote the fitted value in the model space by $\widehat Y$ so that, under the covariance inner product, the residual $Y-\widehat Y$ is orthogonal to (functions in) $\mathcal{H}_{A,X}$.  We then consider the Hilbert space $\mathcal{H}_{X}$ of square integrable functions of $X$, the identical inner product, and consider projecting $A$ --  more generally $S(A)$ from \eqref{eq:EE-S} -- onto $\mathcal{H}_{X}$.  Denote the corresponding fitted value $\widehat A$, and note that the residual $A-\widehat A$ is orthogonal to (functions in) $\mathcal{H}_{X}$.  As $\mathcal{H}_{X} \subseteq \mathcal{H}_{A,X}$, we must then have that $Y-\widehat Y$ is orthogonal to $A-\widehat A$.  This is the fundamental principle that underpins $G$-estimation.

The doubly robust\index{Doubly robust} estimating equation given by \eqref{eq:DREE} then arises naturally from consideration of the orthogonality of the two residuals obtained respectively when fitting $Y$ as a function of $(A,X)$ and $A$ as a function of $X$.

\subsection{A structural version of the outcome model}

The structural, or causal, version of the general model above takes the form
\begin{equation}\label{eq:G-SMM}
\E[Y(a)|\bX] = \E[Y|a,\bX] = m_0(X;\beta) + \gamma(a,X;\theta)
\end{equation}
for treatment $a$, where $Y(a)$ is the counterfactual or potential outcome\index{Potential outcomes} that is typically considered when causal quantities are studied;  here the assumption of no unmeasured confounding,\index{No unmeasured confounding} that is, that $A$ and $Y(a)$ are independent given $X$ is crucial, as it links the potential outcome $Y(a)$ to the observable outcome $Y$ at treatment level $a$ regardless of whether $A$ is regarded as random or fixed by the experimenter.  Causal contrasts are typically based on expected counterfactual outcomes
\[
\E[Y(a)] = \iint y f_{Y|A,X}(y|a,x) f_X(x) \:dy \:dx.
\]
The only difference between this and the standard regression formulation is that treatment arises from an \textit{intervention} rather than as an uncontrolled random quantity.  The inferential complication is that, outside of a randomized trial, data that facilitate estimation of the causal parameters (for example, $\theta$, or contrasts defined explicitly with respect to $\E[Y(a)]$) are not available.  Procedures for estimating these causal quantities rely principally on regression formulations. In the remainder of this paper we retain the counterfactual formulation when considering extensions to the multiple interval case.

\medskip
\noindent \textbf{Terminology:}  The formulations above have essentially led us to consider three model components for the specification of the G-estimating equation, which we have denoted by $m_0(X;\beta)$, $\gamma(a,X;\theta)$ and $\pi_S(X_i;\theta)$.  These models are commonly referred to as the \emph{treatment-free outcome} model, the \emph{blip} model\index{Blip function} and the \emph{treatment} or \emph{propensity score}\index{Propensity score} model.  In the standard G-estimation setting, the blip model is sometimes qualified as the \emph{blip-to-zero} model, to make distinct the idea that it represents an expected (counterfactual) difference in outcome due to treatment relative -- on the additive scale -- to treatment at level zero.

\section{G-estimation for multiple intervals: the structural nested mean model }

The initial use of G-estimation considered semiparametric causal inference\index{Causal inference} for the by-interval effects of a fixed sequence of treatments.  This \emph{standard} application of G-estimation could, in principle, be used to determine an optimal static sequence of treatments under a \emph{myopic} decision strategy\index{Myopic treatment strategy} that seeks to optimize outcome at each successive interval (rather than at the final interval only).  In this section, we describe how standard G-estimation can be applied to perform semiparametric inference.

\subsection{Two intervals}

We consider the counterfactual (terminal) outcome $Y(a_1,a_2)$ arising from the intervention $a_1$ in interval one and $a_2$ in interval two.  We consider the corresponding predictors $\bX_1$ and $\bX_2$, and further denote by $\upto{j}{\bX}$ and $\after{j}{\bX}$ the collection of $\bX$ values \emph{up to and including the $j$th}, and \emph{from the $j$th onward}, respectively; that is $\upto{1}{\bX} = \bX_1$, $\upto{2}{\bX} = (\bX_1,\bX_2)$, $\after{1}{\bX} = (\bX_1,\bX_2)$ and $\after{2}{\bX} = \bX_2$. A key decomposition in understanding the nested mean model expresses $Y(a_1,a_2)$ as follows:
\begin{equation}
Y(a_1,a_2) = Y(a_1,0) + [Y(a_1,a_2) - Y(a_1,0) ]. \label{eq:G-I}
\end{equation}
The first term represents the contribution to the response due to treatment at interval one. The second term represents the incremental response component due to being treated at level $a_2$ at time point two (compared to $a_2 = 0$), after being treated at level $a_1$ at interval one.  Only the sum of these two terms is observed in the data.  Two models, parametrized by $\beta_2$ and $\theta_2$ respectively, are assumed for the conditional expectation of the first term
\begin{equation}
\E[Y(a_1,0)|\upto{2}{\bX}] = m_{02}(\upto{2}{H};\bbeta_2) \label{eq:G-beta2}
\end{equation}
and of the second term
\begin{equation}
\E[Y(a_1,a_2) - Y(a_1,0)|\upto{2}{\bX}] = \gamma_2(a_2,\upto{2}{H};\theta_2) \label{eq:G-gamma2}
\end{equation}
where $\upto{2}{H} = (a_1,\upto{2}{\bX})$, yielding a decomposition similar to that in \eqref{eq:GeneralRegressionForm}.   The parameters $\beta_2,\theta_2$ are estimated using an estimating equation similar to \eqref{eq:EE-S}:
\begin{equation}
\sum_{i=1}^n (Y_i - m_{02}(\upto{i2}{H};\bbeta_2) - \gamma_2(A_{i2},\upto{i2}{H};\theta_2)) \dbinom{T_2(\upto{i2}{H};\bbeta_2)}{S_2(A_{i2},\upto{i2}{H};\theta_2)-\pi_{2S}(\upto{i2}{H};\theta_2,\alpha_2)} = \zerovec
\label{eq:G-II}
\end{equation}
where
\[
T_2(l;\beta) = \frac{\partial m_{02}(l;\beta)}{\partial \beta} \qquad \qquad S_2(a,l;\theta) = \frac{\partial \gamma_2(a,l;\theta)}{\partial \theta}
\]
and
\[
\pi_{2S}(\upto{i2}{H};\theta_2,\alpha_2) = \E[S_2(A_{i2},\upto{i2}{H};\theta_2)|\upto{i2}{H};\alpha]
\]
is the conditional expected value of $S_2(A_{i2},\upto{i2}{H};\theta_2)$ given $\upto{i2}{H}$. Having produced estimates $(\bbetahat_2,\widehat\theta_2)$ for the second interval using \eqref{eq:G-II}, we now proceed to estimate parameters related to the first interval.  From the decomposition \eqref{eq:G-I} we have
\[
Y(a_1,0) = Y(a_1,a_2) - [Y(a_1,a_2) - Y(a_1,0)];
\]
from the fit at the second interval we estimate the second term by $\gamma_2(a_2,\upto{2}{H};\widehat \theta_2)$ for each $i$, and construct the pseudo-outcome
\[
Y_i^{(1)} = \widehat Y(a_1,0) = Y_i - \gamma_2(a_{i2},\upto{i2}{H};\widehat \theta_2).
\]
Then using a similar approach to interval two, we build the regression model for the pseudo-outcome
\begin{equation}
\E[Y_i^{(1)}|a_{i1},\upto{i1}{\bX}] = m_{01}(\upto{i1}{\bX};\bbeta_1) + \gamma_1(A_{i1},\upto{i1}{\bX};\theta_1) \label{eq:G-IIImean}
\end{equation}
and utilize the estimating equation
\begin{equation}
\sum_{i=1}^n (Y_i^{(1)} - m_{01}(\upto{i1}{\bX};\bbeta_1) - \gamma_1(A_{i1},\upto{i1}{\bX};\theta_1)) \dbinom{T_1(\upto{i1}{\bX};\bbeta_1)}{S_1(A_{i1},\upto{i1}{\bX};\theta_1)-\pi_{1S}(\upto{i1}{\bX};\theta_1,\alpha_1)} = \zerovec
\label{eq:G-III}
\end{equation}
to produce estimates $(\bbetahat_1,\widehat \theta_1)$.

\medskip

\noindent \textbf{Note:} We usually proceed and consider a specific form of this general model where the $\gamma(.)$ functions are linear in the parameters $\theta$.  This is the most typically assumed structure, which leads to the random quantities $S_2(A_{i2},\upto{i2}{H};\theta_2)$ and $S_1(A_{i1},\upto{i1}{\bX};\theta_1)$, and their conditional expectations, being functionally independent of $\theta_2$ and $\theta_1$ respectively.

\medskip

\noindent \textbf{Remark:}  The two interval structured model proposes the two related conditional expectations from \eqref{eq:G-beta2} and \eqref{eq:G-IIImean} for the counterfactual quantity $Y(a_1,0)$;  $\E[Y|a_1,\upto{2}{\bX}]$ and $\E[Y|a_{1},\upto{1}{\bX}]$. From first principles we should have
\[
\E[Y|a_{1},\bx_1] = \int \E[Y|a_1,\bx_1,\bx_2] f_{\bX_2|\bX_1} (\bx_2|\bx_1) \: d \bx_2.
\]
This identity implies the need for compatibility between the models
$m_{02}(\upto{2}{H};\bbeta_2)$ and $m_{01}(\bX_1;\bbeta_1) + \gamma_1(a_1,\bX_1;\theta_1)$, that is
\[
\E_{\bX_2|\bX_1}[m_{02}(a,\upto{2}{X};\bbeta_2)] =  m_{01}(\upto{1}{\bX};\bbeta_1) + \gamma_1(a,\upto{1}{\bX};\theta_1).
\]
For example, if we specify
\[
m_{02}(a,\upto{2}{x};\bbeta_2) = \beta_{20} + \beta_{21} x_1 + \beta_{22} x_2 + \beta_{23} a + \beta_{24} a  x_1
\]
then on taking expectations we identify
\[
m_{01}(\upto{1}{\bx};\bbeta_1) = (\beta_{20} + \beta_{22} \xi_{21}(x_1)) + \beta_{21} x_1 \qquad  \gamma_1(a,\upto{1}{\bX};\theta_1) = a(\beta_{23} + \beta_{24} x_1)
\]
where $\xi_{21}(x_1)$ is the conditional expectation of $X_2$ given that $X_1=x_1$.  Thus the specification that suggests that, for example, $m_{01}(\upto{1}{\bx};\bbeta_1)$ is linear in $x_1$ implies that $\xi_{21}(x_1)$ must be either constant or linear in $x_1$.

Furthermore it is clear from extending the decomposition \eqref{eq:G-I} that
\begin{eqnarray*}
Y(a_1,a_2) & = & Y(a_1,0) + [Y(a_1,a_2) - Y(a_1,0) ] \nonumber \\[6pt]
& = & Y(0,0) + [Y(a_1,0) - Y(0,0) ] + [Y(a_1,a_2) - Y(a_1,0)] \label{eq:G-IVa}
\end{eqnarray*}
so the conditional expectations given $\upto{1}{\bX}$ and $\upto{2}{\bX}$ should be compatible, and parameters in a compatible specification are shared between models across intervals.  To guarantee compatibility, we may consider the specifications conditional given $\upto{2}{\bX}$, and marginalization.   Finally, if attention is restricted to a linear model of the form of \eqref{eq:BLM}, we see that the sequential G-estimation strategy is essentially an alternative form of OLS estimation.

\subsection{\texorpdfstring{$K$}{K} intervals}

For the $K$ interval case, we extend the previous procedure to perform inference for the parameters in each interval in reverse order.   Specifically, we first consider the counterfactual outcomes $Y(\upto{K}{a})$ and the decomposition
\begin{eqnarray*}
Y(\upto{K}{a}) &=& Y(\upto{K-1}{a},0)+[Y(\upto{K}{a})-Y(\upto{K-1}{a},0) ] \\[6pt]
&=& Y(\upto{K-2}{a},\zerovec_2) + [Y(\upto{K-1}{a},\zerovec) - Y(\upto{K-2}{a},\zerovec_2)]+[Y(\upto{K}{a})-Y(\upto{K-1}{a},0) ] \\[6pt]
& = & \cdots\\[6pt]
&=& Y(\zerovec_K) + \sum_{l=1}^K \: [Y(\upto{K-l+1}{a},\zerovec_{l-1}) - Y(\upto{K-l}{a},\zerovec_l)]
\end{eqnarray*}
where $\zerovec_k$ is the $k$-dimensional vector of zeros.
%We then construct models
%\[
%\E[Y(\upto{k}{a},\zerovec_{K-k}) - Y(\upto{k-1}{a},\zerovec_{K-k+1}) |\upto{k}{H} ] = \gamma_k (a_k,\upto{k}{H};\theta_k)
%\]
Following the sequential fitting strategy detailed in the previous section, for the $k$th interval, we consider
\[
Y(\upto{k}{a},\zerovec_{K-k}) = Y(\upto{k-1}{a},\zerovec_{K-k+1})+ [Y(\upto{k}{a},\zerovec_{K-k}) -  Y(\upto{k-1}{a},\zerovec_{K-k+1})]
\]
then model the two conditional expectations
\begin{eqnarray*}
\E[Y(\upto{k-1}{a},\zerovec_{K-k+1}) |\upto{k}{H}] &=& m_{0k}(\upto{k}{H};\bbeta_k)\\[6pt]
\E[Y(\upto{k}{a},\zerovec_{K-k}) -  Y(\upto{k-1}{a},\zerovec_{K-k+1})|\upto{k}{H}] &=& \gamma_k (a_k,\upto{k}{H};\theta_k)
\end{eqnarray*}
where $\upto{k}{H} = (\upto{k-1}{a},\upto{k}{\bX})$, and $\gamma_k (0,\upto{k}{H} ;\theta_k) \equiv 0$.  The parameters in these models are estimated using the estimating equation
\begin{equation*}
\sum_{i=1}^n (Y_i^{(k)} - m_{0k}(\upto{ik}{\bX};\bbeta_k) - \gamma_k(A_{ik},\upto{ik}{\bX};\theta_k)) \dbinom{T_k(\upto{ik}{\bX};\bbeta_k)}{S_k(A_{ik},\upto{ik}{\bX};\theta_k)-\pi_{kS}(\upto{ik}{\bX};\theta_k,\alpha_k)} = \zerovec
\label{eq:G-K}
\end{equation*}
where
\begin{equation}\label{eq:psuedok}
Y_i^{(k)} = Y_i - \sum_{l=k+1}^K \gamma_l(a_{il},\upto{il}{H};\widehat\theta_l)
\end{equation}
is the pseudo-outcome for interval $k$.

\section{Optimal treatment regimes via G-estimation}

Our goal is to infer, for each individual in the study, the sequence of treatments that would have (in expectation) led to the optimal outcome for that individual.  We assume that treatment decisions are to be made dynamically, in light of information available -- past treatments, individual characteristics, intermediate responses -- at each time instant.  The inferred sequence of optimal treatments is often referred to as a dynamic treatment regime (DTR)\index{Dynamic treatment regime} for the individual concerned; the inferred optimal DTR will in general vary across individuals.   More generally, the inferred optimal regime can be used to construct optimal personalized regimes for individuals in a future study or trial.

The G-estimation procedures described above can be used to deduce the optimal treatment by interval for individuals in the study.  However, the formulation of one step of the G-estimation procedure needs to be adjusted in the multi-interval setting.  Specifically, the formation of the pseudo-outcome takes into account the inferred optimal treatment in subsequent intervals.  Robins \cite{Robins2004} gives a more general and complete treatment of optimal DTR estimation with SNMMs, estimated by G-estimation. A number of other approaches have been proposed and compared \cite{Murphy2003,MoodieEtAl2007,Rosthoj2006,HendersonEtAl2010}. As a semi-parametric procedure, G-estimation offers several attractive properties, including robustness against some forms of model misspecification \cite{Robins2004}.

In the discrete setting, we suppose that treatment $a$ lies in the finite set of values $\mathcal{A} = \{a_1,\ldots,a_A \}$; for the continuous case, $\calA$ is some bounded real interval. The critical difference between standard G-estimation for SNMMs and G-estimation as used for inferring optimal dynamic regimens is that, although the decomposition in \eqref{eq:G-SMM} is still utilized at each interval, the optimal sequence of treatments must be considered in its entirety.    That is, a treatment that is optimal in a given interval may not be optimal for that interval when the entire sequence of treatments is considered.   We refer to treatments determined as optimal in the context of a given interval only as \emph{myopic}, and those that form part of an optimal `in the long-run' sequence as \emph{hyperopic}.

In earlier sections we have seen how decompositions such as \eqref{eq:G-I} in its general form, and pseudo-outcomes as in \eqref{eq:psuedok}, can be used to estimate interval-specific treatment parameters $\psi_k$. However, for optimal treatment sequences, a different decomposition must be used to construct the pseudo-outcomes.

\subsection{Notation}

As before, we suppose that the observed data consist of $n$ i.i.d. vectors, $O_1, \dots, O_n$, where for each $i$ and over $K$ treatment intervals, $O = (L_0, A_1, \upto{1}{H}, \dots, A_K, \upto{K}{H})$ is a vector of observable quantities. In this vector, $L_0$ denotes a vector of baseline or pre-treatment predictors,
$A_k$ denotes the treatment observed at interval $k$, and
$\upto{k}{H}$ denotes a vector of time-varying predictors associated to interval $k$;  $\upto{k}{H}$ is assumed to be an intermediate response, observed after treatment $A_k$ and prior to treatment $A_{k+1}$, which it may inform, and may contain exogenous or endogenous predictors $\bX_k$ -- that is, $\upto{k}{H} = (\upto{k-1}{A},\upto{k}{X})$.  It is assumed that the treatment $A_{k+1}$ is chosen adaptively according to previous treatments and responses, in particular $\upto{k}{H}$.

We define $Y = g(O)$ to be the final outcome of interest, where $g(\cdot)$ is a suitably defined {\it reward function}, defined such that the larger the value of $Y$, the better the outcome.  Let $\calA_k, k=1, \dots, K$ denote the set of possible treatments at interval $k$, and similarly let $\mathcal{H}_k$ denote the set of possible values for $\upto{k}{H}, k=0, \dots, K$.   A dynamic treatment regime\index{Dynamic treatment regime} $a = (a_1, \dots, a_K)$ is then a sequence of functions $a_k: \upto{k}{H} \to \calA_k$ where
\[
\upto{k}{H} = \mathcal{H}_0
\times \calA_1 \times \mathcal{H}_1 \times \dots \times
\calA_{k-1} \times \mathcal{H}_{k-1}
\]
that map histories to treatments.  We assume the existence of {\it counterfactual} or {\it potential outcomes}\index{Potential outcomes} $Y(a)$ under every possible dynamic regimen $a$, which may be equivalently written $Y(a)$ as the potential treatment sequence $a$ is a deterministic function of $a$. In the DTR context, the potential outcome $Y(a)$ represents the outcome that would have been observed in the event that the treatments $a$ had been assigned according to regimen $a$, rather than the observed regime that led to observed treatment sequence $a$. An optimal hyperopic DTR is one that maximizes the value of the counterfactual outcome in expectation.  Under regularity conditions \cite{Robins2004,Moodie2009_Bka}, the optimal dynamic treatment regime \index{Dynamic treatment regime} is uniquely defined and we write
\[
a^{\opt} = \arg \max_a \E[Y(a)].
\]

\subsection{Hyperopic optimal treatment sequences via dynamic programming}
\index{Dynamic programming|(}
We define the optimal treatment sequence by considering so-called `Bellman optimality' \cite{Bellman}\index{Bellman optimality equation}, and use a dynamic programming approach that proceeds backwards from interval $K$ to interval 1 selecting the optimal treatment at each interval presuming that treatment is \emph{optimally assigned} in subsequent intervals.  That is, we have that $\aopt{\;\;} = (\aopt{1},\ldots,\aopt{K})$ is defined by the realization that
\begin{equation*}\label{eq:Bellman}
\aopt{k} = \arg \max_a \E[Y(a_1,\ldots,a_{k-1},a,\aopt{k+1},\ldots,\aopt{K})].
\end{equation*}
The formulation via expected counterfactual outcomes matches the formulation of earlier sections.  However, the difference here is that for inference we will require data (or pseudo-data) observed under partially optimal treatment sequences; such data are not observed in the study itself, and we need to modify the construction of the pseudo-outcomes compared to \eqref{eq:psuedok}.  In the one interval setting, of course, the optimal hyperopic treatment coincides with the optimal myopic treatment, and no pseudo-data are needed.
\index{Dynamic programming|(}

\subsection{Single interval: discrete treatment options}

In the one interval setting, the optimal treatment is that which maximizes the expected (counterfactual) outcome given the observed data.  In particular, we have for the $i$th individual that
\begin{equation}\label{eq:argmax}
\aopt{i} = \arg \max_{a \: \in \mathcal{A}} \gamma(a,\bX_i;\widehat \psi).
\end{equation}
For example, if treatment is binary so that $\mathcal{A} \equiv \{0,1\}$, with 0 indicating baseline treatment (no treatment or standard care, say), we treat at the non-baseline level if $\gamma(1,\bX_i;\widehat \psi) > \gamma(0,\bX_i;\widehat \psi)$, and do not treat otherwise.   This defines a decision rule for future individuals as a function of the predictors $\bX_i$; the inferred optimal treatment will be different for each individual.   For multi-valued discrete treatment options, we have that the optimal treatment decision is given by \eqref{eq:argmax}; this typically requires an exhaustive search over $\calA$.

\subsection{Two intervals}

In the two interval setting, we revisit the previous decomposition \eqref{eq:G-I}. Recall that
\[
Y(a_1,a_2) = Y(a_1,0) + [Y(a_1,a_2) - Y(a_1,0) ]
\]
so that, to initiate the backward dynamic programming\index{Dynamic programming} strategy, the optimal treatment at interval two is given by
\begin{equation*}\label{eq:argmaxII}
\aopt{i2} = \arg \max_{a \: \in \mathcal{A}_2} \gammaopt{2}(a,\upto{i2}{H};\widehat \psi_2),
\end{equation*}
where $\psi_2$ parameterizes the expected counterfactual difference as in \eqref{eq:G-gamma2} as estimated using \eqref{eq:G-II}.   We then note an alternate decomposition that writes
\begin{eqnarray*}\label{eq:OptSecondInterval}
Y(a_1,\aopt{2}) &=& Y(a_1,a_2) + [Y(a_1,\aopt{2}) - Y(a_1,a_2) ]\\[6pt]
& = & Y(a_1,a_2) + [\{Y(a_1,\aopt{2}) - Y(a_1,0)\} - \{Y(a_1,a_2)-Y(a_1,0)\} ].
\end{eqnarray*}
Inspired by this decomposition, for estimation in the first interval, we consider the pseudo-outcomes
\begin{equation*}\label{eq:PseudoFirst}
\Yhatopt{i}{1} = Y_i + \left[ \gammaopt{2}(\aopt{i2},\upto{i2}{H};\widehat \psi_2) - \gammaopt{2}(a_{i2},\upto{i2}{H};\widehat \psi_2) \right].
\end{equation*}
The first interval pseudo-outcomes are then used to perform inference for the first interval parameters.  We consider the decomposition
\[
Y(a_1,\aopt{2}) = Y(0,\aopt{2}) + [Y(a_1,\aopt{2}) - Y(0,\aopt{2}) ]
\]
and write the expectation of the right hand side similarly to \eqref{eq:G-IIImean} as
\begin{equation}
\E[\Yhatopt{i}{1}|a_{i1},\upto{i1}{\bX}] = \mopt{01}(\upto{i1}{\bX};\bbeta_1) + \gammaopt{1}(a_{i1},\upto{i1}{\bX};\psi_1) \label{eq:G-IIImeanOpt}
\end{equation}
where the superscript `opt' on the functions $m_{01}$ and $\gamma_1$ indicates that these functions assume optimal hyperopic treatment subsequent to interval one.  An estimating equation equivalent to \eqref{eq:G-III} is then used to estimate $\psi_1$, and then the optimal interval one treatment is defined as
\begin{equation*}\label{eq:argmaxI}
\aopt{i1} = \arg \max_{a \: \in \mathcal{A}_1} \gammaopt{1}(a,\upto{i1}{\bX};\widehat \psi_1).
\end{equation*}
yielding the inferred optimal regimen $\aopt{i}= (\aopt{i1},\aopt{i2})$ for each individual $i=1,\ldots,n$, and also the DTR decision function $a$ that defines the optimal allocation rule for future individuals.

\medskip
\noindent \textbf{Terminology:}  The function $\gammaopt{1}(a,X;\theta)$ is also commonly referred to as the \emph{blip} model, but is here usually qualified as the \emph{optimal blip-to-zero treatment} model, to make distinct the idea that it represents an expected (counterfactual) difference in outcome due to treatment relative -- on the additive scale -- to treatment at optimal (rather than zero) level in interval two.\index{Blip function}

\medskip

\noindent \textbf{Note:}  The formulation in \eqref{eq:G-IIImeanOpt} relies upon the assumption that $\gammaopt{1}(a_{i1},\upto{i1}{\bX};\psi_1)$ (and indeed the treatment-free outcome model $\mopt{01}(\upto{i1}{\bX};\bbeta_1)$) is not functionally dependent on $\aopt{2}$: this optimal blip-to-zero formulation assumes that the function $\gammaopt{1}$ takes the same value irrespective of the actual value of $\aopt{2}$.  This assumption limits the scope of the data generating mechanisms that can be correctly represented by this formulation.   However it is necessary in order to facilitate the construction of an implementable optimal dynamic regimen; for new subjects, treatment decisions at interval $k$ can only be made conditional on history up to interval $k$, and cannot depend on future predictors or future optimal treatments, neither of which are available at interval $k$.

\subsection{\texorpdfstring{$K$}{K} intervals}

For the $K$ interval case, we extend the previous procedure to perform inference for the optimal treatment parameters $\psi$ in each interval in reverse order.   For interval $k$, we consider the counterfactual outcome $Y(\upto{k-1}{a},a_k,\afteropt{k+1}{a})$ and two decompositions.  First we write
\begin{equation}
Y(\upto{k-1}{a},a_k,\afteropt{k+1}{a}) = Y(\upto{k-1}{a},0,\afteropt{k+1}{a})+[Y(\upto{k-1}{a},a_k,\afteropt{k+1}{a})-Y(\upto{k-1}{a},0,\afteropt{k+1}{a}) ] \label{eq:recursion1}
\end{equation}
which gives a means of estimation through a model similar to \eqref{eq:G-IIImeanOpt}.  Secondly, we write
\begin{equation}
Y(\upto{k-1}{a},\aopt{k},\afteropt{k+1}{a}) = Y(\upto{k-1}{a},a_k,\afteropt{k+1}{a}) + [Y(\upto{k-1}{a},\aopt{k},\afteropt{k+1}{a})-Y(\upto{k-1}{a},a_k,\afteropt{k+1}{a}) ] \label{eq:recursion2},
\end{equation}
which allows us to construct pseudo-outcomes to be used in the resulting estimating equation.   Specifically, for use in \eqref{eq:recursion1}, we need pseudo versions of $Y(\upto{k-1}{a},a_k,\afteropt{k+1}{a})$.  Using \eqref{eq:recursion2}, we have for $k\leq K-2$
\[
Y(\upto{k-1}{a},a_k,\afteropt{k+1}{a}) = Y(\upto{k+1}{a},\afteropt{k+2}{a}) + [Y(\upto{k-1}{a},a_k,\afteropt{k+1}{a}) - Y(\upto{k+1}{a},\afteropt{k+2}{a}) ]
\]
suggesting the recursion
\[
\Yhatopt{\;\;}{k} = \Yhatopt{\;\;}{k+1} + \left[\gammaopt{k+1}(a_{ik+1}{\hspace*{-1.35em} \raisebox{0.925ex}{\scriptsize{$\textrm{opt}$}}}\;,\upto{ik+1}{H},\widehat \psi_{k+1}) -  \gammaopt{k+1}(a_{ik+1},\upto{ik+1}{H},\widehat \psi_{k+1}) \right]
\]
to define the pseudo-outcome for estimating $\psi_k$. The G-estimation is implemented using the model for the second term in \eqref{eq:recursion1}
\[
\E[Y(\upto{k-1}{a},a,\afteropt{k+1}{a})-Y(\upto{k-1}{a},0,\afteropt{k+2}{a})|\upto{k}{H};\psi_k] = \gammaopt{k}(\upto{k-1}{a},a,\upto{k+1}{H};\psi_k)
\]
as before.   The assumption that $\gammaopt{k}(\upto{k-1}{a},a,\upto{k+1}{H};\psi_k)$ does not depend on the future optimal treatments $\afteropt{k+1}{a}$ is retained here.

To initiate the recursive calculation, we proceed as in the two interval case using the decomposition
\[
Y(\upto{K-1}{a},a_K) = Y(\upto{K-1}{a},0) + [Y(\upto{K-1}{a},a_K) - Y(\upto{K-1}{a},0) ]
\]
and model fitting to estimate parameters $\psi_K$.  We optimize over the possible treatments to obtain $a_K {\hspace*{-.65em} \raisebox{0.925ex}{\scriptsize{$\textrm{opt}$}}}$ where
\[
a_{iK}{\hspace*{-.85em} \raisebox{0.925ex}{\scriptsize{$\textrm{opt}$}}} = \arg \max_{a \: \in \mathcal{A}_K} \gammaopt{k}(a,\upto{iK}{H};\widehat \psi_K),
\]
and form the $(K-1)$st pseudo-outcome as
\[
\Yhatopt{\;\;}{K-1} = Y_i + \left[\gammaopt{K}(a_{iK}{\hspace*{-0.70em} \raisebox{0.925ex}{\scriptsize{$\textrm{opt}$}}}\;,\upto{iK}{H},\widehat \psi_{K}) -  \gammaopt{K}(a_{iK},\upto{iK}{H},\widehat \psi_{K}) \right];
\]
thus the $k$th pseudo-outcome is given by the recursion as
\begin{equation}\label{eq:PseudoYk}
\Yhatopt{i}{k} = Y_i + \sum_{l=k+1}^K
\left[
\gammaopt{l}(a_{il}{\hspace*{-0.70em} \raisebox{0.925ex}{\scriptsize{$\textrm{opt}$}}}\;,\upto{il}{H},\widehat \psi_{l}) -  \gammaopt{l}(a_{il},\upto{il}{H},\widehat \psi_{l}) \right].
\end{equation}

\medskip

\noindent \textbf{Note:}  To reiterate, although the formulation of G-estimation in the standard and DTR settings are closely related, there is a critical difference between the functions $\gamma_k(.;\psi_k)$ and $\gammaopt{k}(.;\psi_k)$.  The former measures the expected difference between counterfactual outcomes due to a treatment $a_k$ compared to absence of treatment, assuming \textbf{absence} of treatment in intervals $k+1,\ldots,K$.  The latter measures the expected difference between counterfactual outcomes due to a treatment $a_k$ compared to absence of treatment, assuming \textbf{optimal} of treatment in \textbf{all} subsequent intervals $k+1,\ldots,K$.  Although we have used identical notation for parameters $\psi_k$, the parameters have different interpretation in the two cases.

\medskip
\noindent \textbf{Terminology:}  The expressions in square brackets in \eqref{eq:PseudoYk}, that is
\[
\mu_k(a, \upto{k}{h};\psi_k) = \gammaopt{k}(\aopt{\;\;}, \upto{k}{h};\psi_k) - \gammaopt{k}(a, \upto{k}{h};\psi_k)
\]
are commonly termed \emph{regret}\index{Regret} functions that measure the penalty in outcome terms for treating suboptimally at level $a$ when the optimal treatment is $\aopt{\;\;}$.

\subsection{Structural nested mean model and G-estimation for a continuous dose}
\label{s:model}

With continuous doses, the need to balance efficacy and tolerability requires more complex modelling than is required for binary doses. If $a$ is a continuous variable then without further constraints, the value of $a$ that maximizes the blip function, the optimal treatment, may potentially be either $+\infty$ or $-\infty$. Thus, the choice of model needs to be reconsidered when treatments are on a continuous scale.  First, one must consider the sensible finite range over which the treatments can realistically be assigned, and bound the treatments during estimation. Second, if it is believed that in some circumstances the optimal dose does not correspond to either extreme of this acceptable range, but rather to some intermediate value, then the functional form of the blip function should not be linear (nor, indeed, monotone) in $a$, but should rather allow for a maximum at an interior point of the interval (\cite{MoodieRichardson2010}). One suitable choice is to assume that the blip function is quadratic in $a$.\index{Blip function} Robins \cite{Robins2004} considered a quadratic blip function model in an
example, while Murphy \cite{Murphy2003} and Rosthoj et al. \cite{Rosthoj2006} also suggest the use
of quadratic models, though in the context of modelling regret functions\index{Regret-regression}; see \cite{MoodieEtAl2007}.

\section{Inferential challenges}

Inference properties for standard G-estimation for structural nested mean models can be obtained from general results for semiparametric inference.  In the one interval setting, the influence function\index{Influence function|see{Influence curve}} corresponding to the G-estimating equation is straightforward to study, and thus the required properties (specifically consistency and asymptotic normality) of the resulting estimators follow under standard assumptions.  The asymptotic variance of the estimators can be expressed via the influence function also, and thus robust (or sandwich) forms of the variance estimator can be derived.  See, for example, \cite{Robins2004} and \cite{vansteelandt2014}.

For G-estimation in the optimal hyperopic DTR setting, the general results from semiparametric inference still dictate the properties of the resulting estimators, but there are several further considerations that are relevant:

\medskip

\noindent \textit{Variance estimation:}  As in the standard setting, doubly robust\index{Doubly robust} versions of the semiparametric estimating procedure require estimation of parameters in the `nuisance' models, that is, the parameters previously listed as $\alpha$ (in the treatment allocation model $\pi(.;.)$) and $\beta$ (in the treatment-free outcome model $\mopt{0j}(.;.)$), and therefore the asymptotic variance for the parameters of interest $\psi$ should must be computed for the adjusted G-estimating equation that uses plug-in versions of the estimated parameters.  See \cite{Robins2004} for technical details of, and \cite{Moodie2009_Bka} for an efficient strategy for computing, the variances from the adjusted estimating equation.

\medskip

\noindent \textit{Non-regularity:}\index{Non-regularity}  Reference \cite{Robins2004} outlines inferential problems that arise from non-regular aspects of the G-estimation procedure in the optimal DTR setting.  Specifically, in the construction of the pseudo-data used in the recursive estimation procedure for interval $k$, there is in general a non-smooth dependence on the parameters $\psi_j$ and their estimators for $j > k$.  For example, in the two interval binary exposure case, if $\aopt{2} = 1 \Longleftrightarrow \psi_{02} + \psi_{12} x_1 > 0$, then there is an optimal decision to treat in the second interval only when the inequality is met.  This decision may choose optimal treatments that are different for the different individuals in the study.  However, if in fact $\psi_{12}$ is close to zero (relative to its standard error), then the optimal treatment decision is subject to instability whenever $x_1$ is in an interval near to the critical ordinate $-\psi_{02}/\psi_{12}$.  The discontinuous nature of the treatment decision, and this instability, renders inference for the $\psi$ parameters in general non-regular; non-standard asymptotic distributions may be obtained, and finite sample behaviour may be highly non-standard.  The problems associated with these `exceptional laws' are described in \cite{Robins2004}, and some possible remedies are outlined in \cite{MoodieRichardson2010} and in \cite{ChakrabortyEtAl2009_SMMR}.

\medskip

\noindent \textit{Model checking:}  Double robustness of the estimating equation \eqref{eq:DREE} appears to be an attractive property as it gives some protection in inference against mis-specification of the treatment or treatment-free outcome models.  However, it has been contended (\cite{Kang2007}) that this double robustness is in practice illusory, as it is usually hard to guarantee in practice that either of the two models is in fact correctly specified, and that deviations from correct specification in both models can lead to serious biases in the estimation of the model parameters, and hence to incorrect identification of the optimal regimens.  It is important therefore to recognize the importance of model checking using, for example, residual plots to verify the orthogonality between outcome and treatment residuals that is a fundamental necessity of the G-estimation procedure -- see equation \eqref{eq:RMN-II}. Residual plotting can also help to verify the critical assumption that the optimal blip-to-zero treatment model $\gammaopt{k}(.;.)$ is correctly specified.  Model checking via residual plotting was described in \cite{Rich2010}.

\section{Illustrations and simulation study}

In our simulations, we consider direct (`forward') simulation from the model using the parameters that determine the contrast functions $\gamma_k(.;.)$, and then determine from the specified models the true values of parameters $\gammaopt{k}(.;.)$.  We parameterize the simulating model using parameter $\lambda$ for clarity.  However, we demonstrate the limitations of this approach, and then construct an alternative simulation strategy based on a different strategy that specifies the blip functions directly.

\subsection{One interval}

\begin{myexample}%\label{ex:Ex1Desc}
We consider first a simple scenario with a single confounder $X_1$ and binary exposure $A_1$, with model specification as follows:
\[
\begin{array}{clcl}
  \textrm{Interval 1} & X_1 & \sim & \Normal(\xi_1,\sigma_1^2) \\[6pt]
    & A_1 |X_1 = x_1 & : & \Pr[A_1 = 1 | X_1 = x_1] = \expit{\alpha_{01} + \alpha_{11} x_1}\\[6pt]
    & Y|A_1=a_1,X_1 = x_1 & \sim &  \Normal(\lambda_{0} + \lambda_1 x_1 + \lambda_2 a_1 ,\sigma_Y^2).
\end{array}
\]
The average potential outcome\index{Potential outcomes} resulting from this structural specification is thus merely
\[
\E[Y(a_1)] = \lambda_{0} + \lambda_1 \xi_1 + \lambda_2 a_1
\]
with average causal effect $\E[Y(1)-Y(0)] = \lambda_2$.  In this model, we see that in the established notation
\begin{eqnarray*}
  m_{01}(a,x_1;\bbeta_1) &=& \beta_{01} + \beta_{11} x_1\\[6pt]
  \gamma_1(a,x_1;\theta_1) &=& \theta_{11} a
\end{eqnarray*}
where $\beta_1 = (\beta_{01},\beta_{11}) \equiv (\lambda_0,\lambda_1)$ and $\theta_1 = (\theta_{11}) = \lambda_2$. In this case it is evident that the optimal treatment $\aopt{1}$ is determined entirely by the sign of $\theta_{11}$ and is identical for each $i$, that is $\aopt{1} = 1 \Longleftrightarrow \theta_{11} > 0$.  Under complete correct specification, G-estimation proceeds using the estimating equation
\[
\sum_{i=1}^n (Y_i - (\beta_{01} + \beta_{11} x_{i1}) - \theta_{11} a_{i1}) \left(\begin{array}{c} 1 \\ x_{i1} \\ a_{i1} - \pi_1(x_{i1};\alpha) \end{array} \right) = \zerovec.
\]
Again, double robustness of the G-estimation procedure ensures that, for example, the estimating equation
\[
\sum_{i=1}^n (Y_i - \beta_{01} - \theta_{11} a_{i1}) \left(\begin{array}{c} 1 \\ a_{i1} - \pi_1(x_{i1};\alpha) \end{array} \right) = \zerovec.
\]
also produces consistent estimation of $\theta_{11}$, as would the estimating equation
\[
\sum_{i=1}^n (Y_i - (\beta_{01} + \beta_{11} x_{i1}) - \theta_{11} a_{i1}) \left(\begin{array}{c} 1 \\ x_{i1} \\ a_{i1} - \pi_1^\ast(x_{i1};\alpha) \end{array} \right) = \zerovec.
\]
for some misspecified treatment model $\pi_1^\ast(x_{i1};\alpha)$.
\end{myexample}

\begin{myexample}%\label{ex:Ex2Desc}
We now consider the model specification identical to the above but with
\[
\begin{array}{cccl}
  \textrm{Interval 1} & Y|A_1=a_1,X_1 = x_1 & \sim &  \Normal(\lambda_{0} + \lambda_1 x_1 + \lambda_2 a_1 + \lambda_3 x_1 a_1 ,\sigma_Y^2)
\end{array}
\]
that contains an interaction between $A_1$ and $X_1$ in the outcome model. The average potential outcome\index{Potential outcomes} resulting from this structural specification is now
\[
\E[Y(a_1)] = \lambda_{0} + \lambda_1 \xi_1 + (\lambda_2 + \lambda_3 \xi_1) a_1
\]
with average causal effect $\E[Y(1)-Y(0)] = \lambda_2 + \lambda_3 \xi_1$.  In this model, we see that in the established notation
\begin{eqnarray*}
  m_{01}(x_1;\bbeta_1) &=& \beta_{01} + \beta_{11} x_1\\[6pt]
  \gamma_1(a,x_1;\theta_1) &=& a(\theta_{11} + \theta_{21} x_1)
\end{eqnarray*}
where $\beta_1 = (\beta_{01},\beta_{11}) \equiv (\lambda_0,\lambda_1)$ and $\theta_1 = (\theta_{11},\theta_{21}) = (\lambda_2,\lambda_3)$. In this case the optimal treatment $\aopt{1}$ is given by
\[
\aopt{1} = 1 \Longleftrightarrow \theta_{11} + \theta_{21} x_1 > 0
\]
which is now potentially different for each individual. Under complete correct specification, G-estimation proceeds using the estimating equation
\[
\sum_{i=1}^n (Y_i - (\beta_{01} + \beta_{11} x_{i1}) - a_{i1}(\theta_{11} + \theta_{21} x_{i1})) \left(\begin{array}{c} 1 \\ x_{i1} \\ a_{i1} - \pi_1(x_{i1};\alpha) \\ x_{i1}(a_{i1} - \pi_1(x_{i1};\alpha)) \end{array} \right) = \zerovec.
\]
Again, double robustness of the G-estimation procedure ensures that misspecification of either the treatment-free outcome model or the exposure model will still result in consistent estimation of the effect of exposure.  For example, if the effect of exposure is modelled using
\[
\gamma_1(a,x_1;\psi_1) = \psi_{11} a
\]
the effect of exposure is still consistently estimated by $\widehat \psi_{11}$.
\end{myexample}

Note that even in the two examples above, replacing the binary exposure by a continuous exposure does not complicate the inferential calculations; the only practical complication is that the search for $\aopt{1}$ requires a search over the support of the exposure distribution.

\subsection{Two intervals}

\begin{myexample}\label{ex:Ex3Desc}
We now consider a scenario with a single confounder in each interval, denoted $X_1$ and $X_2$ respectively, and binary exposure in each interval to yield treatment pattern $(A_1,A_2)$.  We consider model specification as follows:
\[
\begin{array}{clcl}
  \textrm{Interval 1} & X_1 & \sim & \Normal(\xi_1,\sigma_1^2) \\[6pt]
    & A_1 |X_1 = x_1 & : & \Pr[A_1 = 1 | X_1 = x_1] = \expit{\alpha_{01} + \alpha_{11} x_1}\\[12pt]
   \textrm{Interval 2} & X_2|A_1=a_1,X_1 = x_1 & \sim & \Normal(\xi_2(a_1,x_1),\sigma_2^2) \\[6pt]
    & A_2 |A_1=a_1,X_1 = x_1,& : & \Pr[A_2 = 1 | A_1=a_1,X_1 =x_1, X_2 = x_2]\\
     &\qquad X_2 = x_2 && = \expit{\alpha_{02} + \alpha_{12} x_1 + \alpha_{22} x_2 + \alpha_{32} a_1 + \alpha_{42} a_1 x_1 }\\[6pt]
    & Y|A_1=a_1,X_1 = x_1, & \sim &  \Normal(\xi(\upto{2}{a},\upto{2}{x};\lambda) ,\sigma_Y^2)\\[3pt]
    & \qquad A_2=a_2,X_2=x_2 \\
    & & & \xi(\upto{2}{a},\upto{2}{x};\lambda)=\lambda_0 + \lambda_1 a_1 + \lambda_2 a_2 + \lambda_3 x_1 + \lambda_4 x_2\\
\end{array}
\]
In this model, we see that for standard G-estimation beginning at interval 2, a correct outcome mean model has components
\begin{eqnarray*}
  m_{02}(\upto{2}{h};\bbeta_2) &=& \beta_{02} + \beta_{12} a_1 + \beta_{22} x_1 + \beta_{32} x_2\\[3pt]
  \gamma_2(a,\upto{2}{h};\theta_2) &=& \theta_{12} a
\end{eqnarray*}
where $\beta_2 = (\beta_{02},\beta_{12},\beta_{22},\beta_{32}) \equiv (\lambda_0,\lambda_1,\lambda_3,\lambda_4)$ and $\theta_2 = (\theta_{12}) \equiv \lambda_2$.  Again we have that the optimal treatment at interval 2 is determined by the sign of $\theta_{12}$ (or in practice its estimate): $\aopt{2} = 1 \Longleftrightarrow \theta_{12} > 0$.  Note that in interval 2, the correct specification mean model components $\mopt{02}(.;.)$ and $\gammaopt{2}(.;.)$ are identical to $m_{02}(.;.)$ and $\gamma_2(.;.)$ respectively, and
\[
\beta_2 = (\beta_{02},\beta_{12},\beta_{22},\beta_{32}) \equiv (\lambda_0,\lambda_1,\lambda_3,\lambda_4) \qquad \psi_2 \equiv \theta_{12} \equiv \lambda_2.
\]

For interval 1, standard G-estimation proceeds to consider the conditional expectation $\E[Y(a_1,0)|\upto{1}{H}]$ and proposes the decomposition
\[
\E[Y(a_1,0)|\upto{1}{H}] = \E[Y(0,0)|\upto{1}{H}] + \left[\E[Y(a_1,0)|\upto{1}{H}] - \E[Y(0,0)|\upto{1}{H}]\right] .
\]
In the data generating model
\[
\E[Y(a_1,0)|\upto{1}{H}] = \lambda_0 + \lambda_1 a_1 + \lambda_3 x_1 + \lambda_4 \xi_2(a_1,x_1)
\]
where
\[
\xi_2(a_1,x_1) \equiv \E_{X_2|X_1,A_1} [ X_2 | A_1 = a_1, X_1 = x_1]
\]
may involve further parameters that characterize the conditional distribution (or expectation) of $X_2$. Note that in the calculation of $\xi_2(a_1,x_1)$, the relevant expectation is taken with respect to the distribution of $X_2$ given $X_1$, with the value $a_1$ treated as a known fixed constant. Thus for correct specification we require
\begin{eqnarray*}
  m_{01}(\upto{1}{h};\bbeta_1) &=& \beta_{01} + \beta_{11} x_1\\[3pt]
  \gamma_1(a,\upto{1}{h};\theta_1) &=& \theta_{11} a + \theta_{21} \xi_2(a,x_1).
\end{eqnarray*}
Note that the function $\xi_2(a,x_1)$ may be configured such that elements may be passed into function $m_{01}(.;.)$.  Recall also that we require that $\gamma_1(0,.;.) = 0$. Inference then proceeds using the pseudo-outcome
\[
\widehat Y^{(1)}  = Y - \gamma_2(a_{2},\upto{2}{h};\widehat \theta_2),
\]
and an estimating equation specified using functions $m_{01}(.;.)$ and $\gamma_1(.;.)$.

\bigskip

For G-estimation in the optimal DTR setting, the conditional expectation of the outcome assuming optimal treatment in the second interval, $\E[Y(a_1,\aopt{2})|\upto{1}{H}]$, is considered. The decomposition
\[
\E[Y(a_1,\aopt{2})|\upto{1}{H}] = \E[Y(0,\aopt{2})|\upto{1}{H}] + \left[ \E[Y(a_1,\aopt{2})|\upto{1}{H}] - \E[Y(0,\aopt{2})|\upto{1}{H}] \right ]
\]
determines the G-estimating equation, with the first term represented by the function $ \mopt{01}(.;.)$ and the second represented by $\gammaopt{1}(;)$. In the data generating model
\[
\E[Y(a_1,\aopt{2})|\upto{1}{H}= \upto{1}{h}] = \lambda_0 + \lambda_1 a_1 + \lambda_2 \aopt{2} + \lambda_3 x_1 + \lambda_4, \xi_2(a_1,x_1)
\]
meaning that for correct specification we require
\begin{eqnarray*}
  \mopt{01}(\upto{1}{h};\bbeta_1) &=& \beta_{01} + \beta_{11} \aopt{2}+ \beta_{21} x_1\\[3pt]
  \gammaopt{1}(a,\upto{1}{h};\psi_1) &=& \psi_{11} a + \psi_{21} \xi_2(a,x_1).
\end{eqnarray*}
The pseudo-outcome $\Yhatopt{\;\;}{1}$ is then formed using the relation
\[
Y(a_1,\aopt{2}) = Y(a_1,a_2) + [\{Y(a_1,\aopt{2}) - Y(a_1,0)\} - \{Y(a_1,a_2)-Y(a_1,0)\} ],
\]
that is, in estimated form,
\[
\Yhatopt{\;\;}{1} = Y + \left[ \gammaopt{2}(\aopt{2},\upto{2}{H};\widehat \psi_2) - \gammaopt{2}(a_{2},\upto{2}{H};\widehat \psi_2) \right].
\]
Again, the G-estimation procedure is robust to misspecification of these mean models, provided that the expectations of the random quantities
\[
S_k(a,\upto{k}{h};\theta) = \frac{\partial \gammaopt{k}(a,\upto{k}{h};\psi_k)}{\partial \psi_k} \qquad k=1,2
\]
that appear in the G-estimating equations are correctly specified.  In this example,
\[
\beta_1 = (\beta_{01},\beta_{11},\beta_{21}) \equiv (\lambda_0,\lambda_2,\lambda_3) \qquad
\psi_1 = (\psi_{11},\psi_{21}) \equiv (\lambda_1,\lambda_4)
\]
are the correct specification mean model parameters.
\end{myexample}

\medskip
\noindent \textbf{Simulation:} For illustration of Example \ref{ex:Ex3Desc}, we consider the following simulation settings:
\[
\begin{array}{llll}
  \textrm{Interval 1} & X_1 && \xi_1 = 0, \sigma_1 = 1 \\[3pt]
    & A_1|X_1 && \alpha_{01} = -2, \alpha_{11}=2 \\[6pt]
  \textrm{Interval 2} & X_2|A_1,X_1 & \xi_2(a_1,x_1) = \eta_{02} + \eta_{12} a_1 + \eta_{22} x_1 \qquad & \eta_{02}=2,\eta_{12}=-2 , \eta_{22}=1.5,\\
  & && \sigma_2=1 \\[6pt]
  & A_2|A_1,X_1,X_2 && \alpha_{02}=-2,\alpha_{12}=2,\alpha_{22}=1.0,
  \\ & && \alpha_{32}=-1 \\[6pt]
  & Y|A_1 ,X_1  ,A_2 ,X_2  & \multicolumn{2}{l}{\xi(\upto{2}{a},\upto{2}{x};\lambda)=\lambda_0 + \lambda_1 a_1 + \lambda_2 a_2 + \lambda_3 x_1 + \lambda_4 x_2}\\[3pt]
  & & \multicolumn{2}{l}{\qquad \qquad \lambda=(0.2,-1.0,2.0,0.5,0.2), \sigma_Y = 1.}\\
  \end{array}
\]
For this model, the true hyperopic optimal regimen is $(\aopt{1},\aopt{2}) = (0,1)$ for all individuals. For optimal DTR G-estimation, at the second interval, correct model specification sets
\begin{eqnarray*}
  \mopt{02}(\upto{2}{h};\bbeta_2) &=& \beta_{02} + \beta_{12} a_1 + \beta_{22} x_1 + \beta_{32} x_2 \equiv \lambda_0 + \lambda_1 a_1 + \lambda_{3} x_1 + \lambda_{4} x_2\\[3pt]
  \gammaopt{2}(a,\upto{2}{h};\psi_2) &=& \psi_{12} a \equiv \lambda_2 a  .
\end{eqnarray*}
At the first interval, we have that the data generating optimal outcome mean model takes the form
\[
\E[Y(a_1,\aopt{2})|\upto{1}{H}= \upto{1}{h}] = \lambda_0 + \lambda_1 a_1 + \lambda_2 \aopt{2} + \lambda_3 x_1 + \lambda_4 (\eta_{02} + \eta_{12} a_1 + \eta_{22} x_1)
\]
and thus after collecting terms we have
\begin{eqnarray*}
  \mopt{01}(\upto{1}{h};\bbeta_1) &=& \beta_{01} + \beta_{11} \aopt{2}+ \beta_{21} x_1\\[3pt]
  \gammaopt{1}(a,\upto{1}{h};\psi_1) &=& \psi_{11} a
\end{eqnarray*}
where
\[
\beta_{01} = \lambda_0 + \lambda_4 \eta_{02} \qquad \beta_{11} = \lambda_2 \qquad \beta_{21} = \lambda_3 + \lambda_4 \eta_{22}
\]
and
\[
\psi_{11} = \lambda_1 + \lambda_4 \eta_{12}.
\]
In this example, as $\aopt{2}$ is identical for all individuals in the data generating process and so after inference at the second interval the parameters $(\beta_{01},\beta_{11})$ in $\mopt{01}(.;.)$ are not individually identifiable.

For illustration of the G-estimation procedure, we fit misspecified versions of these outcome models, but correctly specify the parametric treatment models $\pi_{1}(.;.)$ and $\pi_{2}(.;.)$.  We consider correct and incorrect specification of the treatment-free outcome model components, and also consider cases where the optimal blip-to-zero\index{Blip function} treatment model is over-specified (that is, where the data generating model is nested within the fitted model).

\bigskip

\textbf{Fitted Models:}
\begin{itemize}
  \item Fitted model 1:  Correct specification.
  \[
  \begin{array}{llcll}
    \mopt{02} & =  \beta_{02} + \beta_{12} a_1 + \beta_{22} x_1 + \beta_{32} x_2 & \qquad \qquad & \mopt{01} &  = \beta_{01} + \beta_{11} x_1 \\[6pt]
    \gammaopt{2} & =  \psi_{12} a_2  & \qquad \qquad & \gammaopt{1} & = \psi_{11} a_1
   \end{array}
  \]

  \item Fitted model 2:  Incorrect specification of treatment-free outcome model in Interval 2.
  \[
  \begin{array}{llcll}
    \mopt{02} & =  \beta_{02} + \beta_{12} a_1 + \beta_{22} x_1 & \qquad \qquad & \mopt{01} &  = \beta_{01} + \beta_{11} x_1 \\[6pt]
    \gammaopt{2} & =  \psi_{12} a_2  & & \gammaopt{1} & = \psi_{11} a_1
  \end{array}
  \]

\item Fitted model 3:  Incorrect specification of treatment-free outcome model in both intervals.
  \[
  \begin{array}{llcll}
    \mopt{02} & =  \beta_{02} + \beta_{12} a_1 + \beta_{22} x_1 & \qquad \qquad & \mopt{01} &  = \beta_{01} \\[6pt]
    \gammaopt{2} & =  \psi_{12} a_2  & \qquad \qquad & \gammaopt{1} & = \psi_{11} a_1
  \end{array}
  \]

  \item Fitted model 4:  Incorrect specification of optimal blip-to-zero treatment outcome model in both intervals
  \[
  \begin{array}{llcll}
    \mopt{02} & =  \beta_{02} + \beta_{12} a_1 + \beta_{22} x_1 & \qquad \qquad & \mopt{01} &  = \beta_{01} + \beta_{11} x_1 \\[6pt]
    \gammaopt{2} & =  \psi_{12} a_2 + \psi_{22} a_2 x_2 & \qquad \qquad &
    \gammaopt{1} & = \psi_{11} a_1 + \psi_{21} a_1 x_1
  \end{array}
  \]

  \item Fitted model 5:  Incorrect specification of both models in both intervals
  \[
  \begin{array}{llcll}
    \mopt{02} & =  \beta_{02} + \beta_{12} a_1 + \beta_{22} x_1 & \qquad \qquad & \mopt{01} &  = \beta_{01}\\[6pt]
    \gammaopt{2} & =  \psi_{12} a_2 + \psi_{22} a_2 x_2 & \qquad \qquad & \gammaopt{1} & = \psi_{11} a_1 + \psi_{21} a_1 x_1
  \end{array}
  \]

  \item Fitted model 6:  Incorrect specification of both models in both intervals
  \[
  \begin{array}{llcll}
    \mopt{02} & =  \beta_{02} & \qquad \qquad & \mopt{01} &  = \beta_{01}\\[6pt]
    \gammaopt{2} & =  \psi_{12} a_2 + \psi_{22} a_2 x_2 & \qquad \qquad & \gammaopt{1} & = \psi_{11} a_1 + \psi_{21} a_1 x_1
  \end{array}
  \]
\end{itemize}

Table \ref{tab:Ex3BiasVar} contains the results of a Monte Carlo study: 5000 replicate data sets were generated using the parameter specification above with $n=200$: results represent the scaled bias and variance for estimators from the six models.  It is clear that all estimators are unbiased, demonstrating the robustness to mis-specification of the outcome model, but that the variance of estimator is influenced by the degree of mis-specification.

\begin{table}[ht]
\centering
\caption{Example 3: Results from 5000 replicate data sets of size $n=200$. }
\label{tab:Ex3BiasVar}
\begin{tabular}{|r|c|r|r|r|r|r|r|r|}
  \hline
&&& \multicolumn{6}{c|}{Model}\\
&Parameter    & \multicolumn{1}{c|}{Truth} & \multicolumn{1}{c|}{1}  &  \multicolumn{1}{c|}{2} & \multicolumn{1}{c|}{3}  & \multicolumn{1}{c|}{4}& \multicolumn{1}{c|}{5}  & \multicolumn{1}{c|}{6}  \\
  \hline
$\sqrt{n} \times$Bias &$\psi_{12}$ & 2.0  & -0.015 & -0.014 &  -0.014 & 0.117  & 0.779   & 2.478  \\
&$\psi_{22}$ & 0.0  & --  & --  & --  & -0.145  & -1.017   &  -1.071 \\[6pt]
&$\psi_{11}$ & -1.4 & 0.075 &  0.075 & -1.017  & 0.135   & 0.135 &  1.505 \\
&$\psi_{21}$ & 0.0  & --  & --  & --  & -0.042  & -0.042  & -1.084  \\
 \hline
$n \times$Var. &$\psi_{12}$ & 2.0  & 11.751 &  11.752 &  11.752 & 42.113  & 47.755  & 71.340  \\
&$\psi_{22}$ & 0.0  & --  & --  & --  & 27.699  & 50.073  &  47.448 \\[6pt]
&$\psi_{11}$ & -1.4 & 16.954&  16.958 & 16.958  & 34.949  & 34.949 &  57.659 \\
&$\psi_{21}$ & 0.0  & --  & --  & --  & 9.918  & 9.918  & 21.634  \\
  \hline
\end{tabular}
\vspace{0.1in}
\end{table}

\begin{myexample} \label{ex:Ex4Desc}
If the model in Example 3 is modified to
\[
\begin{array}{cccl}
  \textrm{Interval 2} & Y|\upto{2}{A}=\upto{2}{a},\upto{2}{X}=\upto{2}{x} & \sim &  \Normal(\xi(\upto{2}{a},\upto{2}{x};\lambda) ,\sigma_Y^2)\\[3pt]
\end{array}
\]
with
\begin{equation}\label{eq:Ex4model}
\xi(\upto{2}{a},\upto{2}{x};\lambda)=\lambda_0 + \lambda_1 a_1 + \lambda_2 a_2 + \lambda_3 x_1 + \lambda_4 x_2 + \lambda_5 a_1 x_1 + \lambda_6 a_2 x_2
\end{equation}
then we have at interval 2 that
\begin{eqnarray*}
  m_{02}(\upto{2}{h};\bbeta_2) &=& \beta_{02} + \beta_{12} a_1 + \beta_{22} x_1 + \beta_{32} x_2 + \beta_{42} a_1 x_1 \\[3pt]
  \gamma_2(a,\upto{2}{h};\theta_2) &=& a (\theta_{12} + \theta_{22} x_2)
\end{eqnarray*}
where $\beta_2 = (\beta_{02},\beta_{12},\beta_{22},\beta_{32},\beta_{42}) \equiv (\lambda_0,\lambda_1,\lambda_3,\lambda_4,\lambda_5)$ and $\theta_2 = (\theta_{12},\theta_{22}) \equiv (\lambda_2,\lambda_6)$.  Then $\aopt{2} = 1 \Longleftrightarrow \theta_{12} + \theta_{22} x_2 > 0.$

For interval 1, we have for standard G-estimation that the correct model specification sets
\begin{eqnarray*}
  m_{01}(\upto{1}{h};\bbeta_1) &=& \beta_{01} + \beta_{11} x_1\\[3pt]
  \gamma_1(a,\upto{1}{h};\theta_1) &=& a(\theta_{11} + \theta_{21} \xi_2(a,x_1))
\end{eqnarray*}
where $\beta_1 = (\beta_{01},\beta_{11}) \equiv (\lambda_0,\lambda_3)$ and $\theta_1 = (\theta_{11},\theta_{21}) \equiv (\lambda_1,\lambda_5)$.  For the optimal DTR version of G-estimation, we have that
\begin{eqnarray*}
  \mopt{01}(\upto{1}{h};\bbeta_1) &=& \beta_{01} +  \beta_{11} \aopt{2} + \beta_{21} x_1  \\[3pt]
  \gammaopt{1}(a,\upto{1}{h};\psi_1) &=& a(\psi_{11} + \psi_{21} x_1) + \psi_{31} \xi_2(a,x_1) + \psi_{41} \aopt{2} \xi_2(a,x_1)
\end{eqnarray*}
where $\beta_1 = (\beta_{01},\beta_{11},\beta_{21}) \equiv (\lambda_0,\lambda_2,\lambda_3)$ and $\psi_1 = (\psi_{11},\psi_{21},\psi_{31},\psi_{41}) \equiv (\lambda_1,\lambda_5,\lambda_4,\lambda_6)$.  Depending on the form of $\xi_2(a,x_1)$, it may be that terms can be further separated.  Note also that, as the optimal blip-to-zero\index{Blip function} function cannot depend on $\aopt{2}$, we must have either that $\psi_{41} \equiv \lambda_6 = 0$, or that $\xi_2(a,x_1)$ does not depend on $a$, or both.  The former assumption corresponds to a lack of interaction between $A_2$ and $X_2$ in the original conditional model \eqref{eq:Ex4model}, which has implications for the second interval analysis above, namely that there can be no tailoring of the second interval treatment by level of $X_2$.  The latter assumption implies that $A_1$ does not influence $X_2$.  Therefore, neither assumption seems easy to justify for practical purposes.
\end{myexample}

\medskip
\noindent \textbf{Simulation:} For Example \eqref{ex:Ex4Desc}, we consider identical simulation settings apart from the final conditional expectation model, and the model relating first interval exposure to second interval confounders.  Specifically, we set
\[
\begin{array}{llll}
  \textrm{Interval 1} & X_1 && \xi_1 = 0, \sigma_1 = 1 \\[3pt]
    & A_1|X_1 && \alpha_{01} = -2, \alpha_{11}=2 \\[6pt]
  \textrm{Interval 2} & X_2|A_1,X_1 & \xi_2(a_1,x_1) = \eta_{02} + \eta_{12} a_1 + \eta_{22} x_1 \qquad & \eta_{02}=2,\eta_{12}=0 , \\
  &&& \eta_{22}=1.5,\sigma_2=1 \\[6pt]
  & A_2|A_1,X_1,X_2 && \alpha_{02}=-2,\alpha_{12}=2,\\
  &&& \alpha_{22}=1.0,\alpha_{32}=-1 \\[6pt]
  & Y|A_1 ,X_1  ,A_2 ,X_2  & \multicolumn{2}{l}{\xi(\upto{2}{a},\upto{2}{x};\lambda)=\lambda_0 + \lambda_1 a_1 + \lambda_2 a_2 + \lambda_3 x_1 + \lambda_4 x_2 + \lambda_5 a_1 x_1 + \lambda_6 a_2 x_2 }\\[3pt]
  & & \multicolumn{2}{l}{\qquad \qquad \lambda=(0.2,-1.0,2.0,0.5,0.2,1.5,-0.5), \sigma_Y = 1.}\\
  \end{array}
\]
For this model, the true hyperopic optimal regimen is dependent on the values of the predictors. At the second interval, correct model specification would set
\begin{eqnarray*}
  \mopt{02}(\upto{2}{h};\bbeta_2) &=& \beta_{02} + \beta_{12} a_1 + \beta_{22} x_1 + \beta_{32} x_2 +
  \beta_{42} a_1 x_1 \equiv \lambda_0 + \lambda_1 a_1 + \lambda_{3} x_1 + \lambda_{4} x_2 + \lambda_5 a_1 x_1\\[3pt]
  \gammaopt{2}(a,\upto{2}{h};\psi_2) &=& a (\psi_{12} + \psi_{22} x_2) \equiv a (\lambda_2  +  \lambda_6 x_2) .
\end{eqnarray*}
At the first interval, the data generating optimal outcome mean model takes the form
\begin{eqnarray*}
\E[Y(a_1,\aopt{2})|\upto{1}{H}= \upto{1}{h}] & = & \lambda_0 + \lambda_1 a_1 + \lambda_2 \aopt{2} + \lambda_3 x_1 + \lambda_4 (\eta_{02} + \eta_{12} a_1 + \eta_{22} x_1)\\[6pt]
&& + \lambda_5 a_1 x_1 + \lambda_6 \aopt{2} (\eta_{02} + \eta_{12} a_1 + \eta_{22} x_1)
\end{eqnarray*}
and thus after collecting terms we have
\begin{eqnarray*}
  \mopt{01}(\upto{1}{h};\bbeta_1) &=& \beta_{01} + \beta_{11} \aopt{2}+ \beta_{21} x_1 + \beta_{31} \aopt{2} x_1 \\[3pt]
  \gammaopt{1}(a,\upto{1}{h};\psi_1) &=& a (\psi_{11} + \psi_{21} x_1 + \psi_{31} \aopt{2})
\end{eqnarray*}
where
\[
\beta_{01} = \lambda_0 + \lambda_4 \eta_{02} \qquad \beta_{11} = \lambda_2 + \lambda_6
 \eta_{02} \qquad \beta_{21} = \lambda_3 + \lambda_4 \eta_{22} \qquad \beta_{31} = \lambda_6
 \eta_{22}
\]
and
\[
\psi_{11} = \lambda_1 + \lambda_4 \eta_{12} \qquad \psi_{21} = \lambda_5 \qquad \psi_{31} = \lambda_6 \eta_{12}.
\]
With $\eta_{12}$ set equal to zero, the optimal blip-to-zero treatment model $\gammaopt{1}(a,\upto{1}{h};\psi_1)$ does not depend on $\aopt{2}$, as required.

In this example, we study first the accuracy of determination of the optimal regime from a correctly specified optimal blip-to-zero model (that is, for a given predictor profile, the proportion of times that the estimated optimal regime coincides with the true optimal regime), as well as the variance for the estimators of the parameters $\psi_1 = (\psi_{11},\psi_{21})$ and $\psi_2 = (\psi_{12},\psi_{22})$.  Correct optimal blip-to-zero and treatment models were fitted, with mis-specified versions of the treatment-free outcome models:
\[
\mopt{02}  :  \quad \beta_{02} + \beta_{12} a_1 + \beta_{22} x_1 \qquad \qquad \mopt{01}   : \quad \beta_{01} + \beta_{11} x_1.
\]

\begin{table}[ht]
  \centering
\caption{Example 4: Proportion of accurate identification of the true optimal regime, and $n  \times \textrm{Variance}$ of G-estimators computed using 5000 replicate data sets and five sample sizes under correct model specification}\label{tab:Ex4a}
 \begin{tabular}{|lr|r|r|r|r|r|r|}
    \hline
&&\multicolumn{5}{c|}{Sample size $n$}\\[3pt]
&&100&200&500&1000&2000\\
\hline
\multicolumn{2}{|r|}{Proportion accuracy}& 0.750&0.864&0.916&0.941&0.958\\[3pt]
\hline
$n \times \textrm{Var.}$ & $\psi_{11}$ &192.56&40.508&17.81&17.415&17.367\\
&$\psi_{21}$&316.348&56.150&23.575&23.010&22.980\\
&$\psi_{12}$&10369.640&3586.284&74.294&67.478&64.546\\
&$\psi_{22}$&2595.538&1234.684&9.889&8.571&8.232\\
    \hline
  \end{tabular}
  \vspace{0.1 in}
\end{table}
For these parameter settings, the procedure accurately identifies the true optimal regime with reasonably high probability.

For a second analysis, we now consider fitting a mis-specified exposure model at interval 2, allowing for dependence only on $X_2$ and not $(A_1,X_1,X_2)$ as in a correctly specified model.   As the treatment-free outcome model $\mopt{02}(.;.)$ is also mis-specified, the resulting estimators are no longer consistent.  Table \ref{tab:Ex4b} contains scaled bias and variance for the resulting estimators.
\begin{table}[ht]
\caption{Example 4: Proportion accuracy for identifying the true optimal regime, and $\sqrt{n}  \times \textrm{Bias}$ of G-estimators computed using 5000 replicate data sets and five sample sizes under a mis-specified exposure model.}\label{tab:Ex4b}
  \centering
 \begin{tabular}{|lr|r|r|r|r|r|r|}
    \hline
&&\multicolumn{5}{c|}{Sample size $n$}\\[3pt]
&&100&200&500&1000&2000\\
\hline
\multicolumn{2}{|r|}{Proportion accuracy}& 0.133&0.127&0.136&0.141&0.143\\[3pt]
\hline
$\sqrt{n} \times \textrm{Bias}$ & $\psi_{11}$ &-8.766&-12.942&-20.973&-30.195&-43.249\\
&$\psi_{21}$&13.538&19.396&31.515&45.338&64.880\\
&$\psi_{12}$&25.027&21.096&32.303&46.355&66.009\\
&$\psi_{22}$&-6.249&-4.623&-6.701&-9.482&-13.405\\
\hline
$n \times \textrm{Var.}$ & $\psi_{11}$ &289.559&24.977&17.160&16.666&16.628\\
&$\psi_{21}$&643.068&30.606&22.497&22.044&21.981\\
&$\psi_{12}$&8964.003&471.148&93.100&86.864&81.049\\
&$\psi_{22}$&2353.984&73.202&13.553&12.207&11.711\\
    \hline
  \end{tabular}
  \vspace{0.1 in}
\end{table}
The correct regime is identified with low probability, even when the sample size is large.

\medskip

\subsection{Regret-based formulation of the data generating model}
\label{sec:Regretbased}\index{Regret|(}
The interpretation of $\gammaopt{k}(.;\psi_k)$ as an expectation of a counterfactual outcome contrast assuming optimal treatment after interval $k$ makes simulation of outcomes from the model for desired $\psi_k$s more complicated than for the standard structural mean model.  That is, the optimal blip-to-zero\index{Blip function} parameters can be deduced from a full conditional specification, provided the specification obeys the rules that allow, for example, $\gammaopt{k}(.;\psi_k)$ to be correctly specified and not depend on future values, but it is not as straightforward to construct a conditional model suitable for data generation for which the induced blip parameters equal the desired values.  One strategy (see \cite{Murphy2003}) involves direct simulation of the counterfactual outcome under an optimal treatment sequence, and the computation of the corresponding outcome for any other treatment sequence using the recursion given in relation in \eqref{eq:recursion1}, which depends explicitly on the blip parameters.   For example, for two intervals, we have
\begin{eqnarray*}
Y(\aopt{1},\aopt{2}) &=&  Y(a_1,\aopt{2}) + \left[Y(\aopt{1},\aopt{2}) - Y(a_1,\aopt{2}) \right]  \\
   &=& Y(a_1,\aopt{2}) + \left[\left\{ Y(\aopt{1},\aopt{2}) - Y(0,\aopt{2}) \right\} - \left\{ Y(a_1,\aopt{2}) - Y(0,\aopt{2}) \right\} \right]\\[6pt]
\Longleftrightarrow  Y(a_1,\aopt{2}) &=& Y(\aopt{1},\aopt{2}) - \left[\left\{ Y(\aopt{1},\aopt{2}) - Y(0,\aopt{2}) \right\} - \left\{ Y(a_1,\aopt{2}) - Y(0,\aopt{2}) \right\} \right].
\end{eqnarray*}
Under the chosen model specification, the term in square brackets has conditional expectation
\[
\gammaopt{1}(\aopt{1},\upto{1}{h};\psi_1) - \gammaopt{1}(a_1,\upto{1}{h};\psi_1).
\]
Once the optimal regimen $(\aopt{1},\aopt{2})$ is defined, if a simulated variate for $Y(\aopt{1},\aopt{2})$ is available and is denoted $y(\aopt{1},\aopt{2})$, then we may generate a value for $Y(a_1,\aopt{2})$ for a specified $a_1$ as
\begin{eqnarray*}
y(a_1,\aopt{2}) &=& y(\aopt{1},\aopt{2}) - \left[\gammaopt{1}(\aopt{1},\upto{1}{h};\psi_1) - \gammaopt{1}(a_1,\upto{1}{h};\psi_1) \right]\\[6pt]
& \equiv & \gammaopt{1}(a_1,\upto{1}{h};\psi_1) + \left(y(\aopt{1},\aopt{2}) - \gammaopt{1}(\aopt{1},\upto{1}{h};\psi_1)  \right).
\end{eqnarray*}
We then also have that
\begin{eqnarray*}
Y(a_1,\aopt{2}) &=&  Y(a_1,a_2) + \left[Y(a_1,\aopt{2}) - Y(a_1,a_2) \right]  \\
   &=& Y(a_1,a_2) + \left[\left\{ Y(a_1,\aopt{2}) - Y(a_1,0) \right\} - \left\{ Y(a_1,a_2) - Y(a_1,0) \right\} \right]\\[6pt]
\Longleftrightarrow  Y(a_1,a_2) &=& Y(a_1,\aopt{2}) - \left[\left\{ Y(a_1,\aopt{2}) - Y(a_1,0) \right\} - \left\{ Y(a_1,a_2) - Y(a_1,0) \right\} \right]
\end{eqnarray*}
so that, using the previous logic, we may finally obtain a variate $y(a_1,a_2)$ as
\begin{eqnarray*}
y(a_1,a_2) &=& y(a_1,\aopt{2}) - \left[ \gammaopt{2}(\aopt{2},\upto{2}{h};\psi_2) - \gammaopt{2}(a_2,\upto{2}{h};\psi_2) \right]\\[6pt]
& = & y(\aopt{1},\aopt{2}) - \left[\gammaopt{1}(\aopt{1},\upto{1}{h};\psi_1) - \gammaopt{1}(a_1,\upto{1}{h};\psi_1) \right] \\[6pt]
&& \qquad \qquad \qquad - \left[ \gammaopt{2}(\aopt{2},\upto{2}{h};\psi_2) - \gammaopt{2}(a_2,\upto{2}{h};\psi_2) \right].
\end{eqnarray*}
In the $K$ interval setting, we may extend this to the form
\begin{equation}\label{eq:MurphySim}
y(\upto{K}{a}) = y(\uptoopt{K}{a}) - \sum_{k=1}^K \left[ \gammaopt{k}(\aopt{k},\upto{k}{h};\psi_k) - \gammaopt{k}(a_k,\upto{k}{h};\psi_k) \right].
\end{equation}
Thus to simulate a random outcome in accordance with a specified collection of optimal blip-to-zero functions and parameters, we proceed as follows:
\begin{enumerate}

\item Simulate the pattern of treatment and confounder variables for intervals $1,2,\ldots,K$, $(\upto{K}{a},\upto{K}{X})$, according to a specified stochastic model.

\item Directly simulate the outcome under optimal treatment to obtain $y(\uptoopt{K}{a})$.  Typically, this is performed conditional on baseline predictors $X_1$ only.

\item At each interval $k=1,2,\ldots,K$, compute the optimal hyperopic treatment $\aopt{k}$ as
\[
\aopt{k} = \arg \max_a \gammaopt{k}(a,\upto{k}{h};\psi_k)
\]
based on the history up to interval $k$ and the desired parameters $\psi_k$.

\item For the simulated treatment pattern $\upto{K}{a}$ from Step 1, compute $y(\upto{K}{a})$ from $y(\uptoopt{K}{a})$ according to \eqref{eq:MurphySim}.

\end{enumerate}
\noindent In this formulation, the treatment-free outcome models across the intervals are not specified explicitly but may be deduced.

\medskip
\begin{myexample}\label{ex:Ex5desc}
For $K=2$ with binary treatment, suppose for illustration that we require that
\begin{eqnarray*}
\gammaopt{1}(a,\upto{1}{h};\psi_1) & = & a (\psi_{11} + \psi_{21} x_1 + \psi_{31} x_1^2) =  a \bX_1 \psi_1, \\[6pt]
\gammaopt{2}(a,\upto{2}{h};\psi_2) & = & a (\psi_{12} + \psi_{22} x_1 + \psi_{32} x_1 a_1 + \psi_{42} x_2 + \psi_{52} x_1 x_2) = a \bX_2 \psi_2
\end{eqnarray*}
say, where $\bX_1 = (1,x_1,x_1^2)$ and $\bX_2 = (1,x_1,x_1 a_1, x_2,x_1 x_2)$, so that treatment is tailored by $X_1$ in interval 1, and by $(X_1,A_1,X_2)$ in interval 2.  Suppose also that we specify
\[
\E[Y(\aopt{1},\aopt{2})|X_1=x_1] = \lambda_0 + \lambda_1 x_1 + \lambda_2 x_1^2.
\]
Under the same data generating mechanism for the exposure/treatment variables as in Example \ref{ex:Ex3Desc}
%\[
%\begin{array}{llcl}
%  \textrm{Interval 1} & X_1 & \sim & \Normal(\xi_1,\sigma_1^2), \; \xi_1 = %0, \sigma_1 = 1 \\[6pt]
%    & A_1 |X_1 = x_1 & : & \Pr[A_1 = 1 | X_1 = x_1] = \expit{\alpha_{01} + %\alpha_{11} x_1}\\[3pt]
%    &&& \alpha_{01} = -2, \alpha_{11}=2 \\[12pt]
%   \textrm{Interval 2} & X_2|A_1=a_1,X_1 = x_1 & \sim & %\Normal(\xi_2(a_1,x_1),\sigma_2^2) \\[6pt]
%   &&& \xi_2(a_1,x_1) = \eta_{02} + \eta_{12} a_1 + \eta_{22} x_1\\[3pt]
%   &&& \eta_{02}=2,\eta_{12}=-2 , \eta_{22}=1.5,\sigma_2=1 \\[6pt]
%    & A_2 |A_1=a_1,X_1 = x_1, X_2 = x_2 & : & \Pr[A_2 = 1 | A_1=a_1,X_1 %=x_1, X_2 = x_2]\\[3pt]
%     &&& = \expit{\alpha_{02} + \alpha_{12} x_1 + \alpha_{22} x_2 + %\alpha_{32} a_1 + \alpha_{42} a_1 x_1 }\\[3pt]
%&&&   \alpha_{02}=-2,\alpha_{12}=2,\alpha_{22}=1.0,\alpha_{32}=-1
%\end{array}
%\]
we may simulate outcomes readily as in the code in the listing below.

\begin{lstlisting}[title=\textbf{Code for regret-based simulation},label=list:TwoInt]

n<-1000  #Sample size
mu1<-0;sig1<-1;X1<-rnorm(n,mu1,sig1)  #For X1

al01<--2.0;al11<-2.0
A1<-rbinom(n,1,expit(al01+al11*X1))  #For A1

eta02<-2.0;eta12<--2.0;eta22<-1.5;sig2<-1
X2<-rnorm(n,eta02+eta12*A1+eta22*X1,sig2)  #For X2

al02<--2.0;al12<-2.0;al22<-1.0;al32<--1.0
A2<-rbinom(n,1,expit(al02+al12*X1+al22*A1+al32*X2)) #For A2

#Simulate Yopt
lam<-c(0.2,-1.0,2.0)                   #Parameters
Xm<-cbind(1,X1,X1^2)                   #Design matrix
muY.opt<-Xm %*%lam;sigY<-1.0           #Linear predictor
Yopt<-rnorm(n,muY.opt,sigY)

#First stage optimal blip-to-zero
psi.1<-c(1.0,2.0,-0.5)                 #Parameters
Xm1<-cbind(1,X1,X1^2)                  #Design matrix
A1.opt.act<-rep(1,n)*(Xm1 %*% psi.1>0) #Actual optimal treatment

#Second stage optimal blip-to-zero
psi.2<-c(1.0,0.2,0.5,-0.4,-0.2)        #Parameters
Xm2<-cbind(1,X1,X1*A1,X2,X1*X2)        #Design matrix
A2.opt.act<-rep(1,n)*(Xm2 %*% psi.2>0) #Actual optimal treatment

#Simulated outcome
Y<-Yopt-(Xm1%*%psi.1)*(A1.opt.act-A1)-(Xm2%*%psi.2)*(A2.opt.act-A2)
\end{lstlisting}
In this data generating model, it is evident that the distribution of the outcome conditional on $(X_1,A_1,X_2,A_2)$ is normal with mean
\begin{equation}\label{eq:MurphySim2}
(\lambda_0 + \lambda_1 x_1 + \lambda_2 x_1^2) - \bX_1 \psi_1 (\aopt{1}-a_1) - \bX_2 \psi_2 (\aopt{2}-a_2)
\end{equation}
where
\[
\aopt{1} = 1_{\{\bX_1 \psi_1 > 0 \}} \qquad \aopt{2} = 1_{\{\bX_2 \psi_2 > 0 \}}.
\]
It is clear, however, that the non-linear dependence of \eqref{eq:MurphySim2} on $X_1$ and $X_2$ renders the calculation of the implied treatment free outcome models more complicated than in the previous formulation.  Therefore, correct specification for these models is unlikely in practice.
\end{myexample}

\medskip

\noindent \textbf{Note:}  Suppose that we study the more general data generating procedure that allows for
\[
\E[Y(\aopt{1},\aopt{2})|\upto{2}{X}] = \lambda_0 + \lambda_1 x_1 + \lambda_2 x_2
\]
say, that is, dependent on $X_1$ and $X_2$.  Then estimation at interval 1 is based on the conditional expectation decomposition
\begin{eqnarray*}
  \E[Y(a_1,\aopt{2})|\upto{1}{l}] &=& \E[Y(\aopt{1},\aopt{2})|\upto{1}{l}] - [ \gammaopt{1}(\aopt{1},\upto{1}{l};\psi_1) - \gammaopt{1}(a_1 ,\upto{1}{l};\psi_1)] \\[6pt]
   &=& \E[(\lambda_0 + \lambda_1 X_1 + \lambda_2 X_2)|\upto{1}{l}] - \bX_1 \psi_1 (\aopt{1}-a_1)\\[6pt]
   &=& \lambda_0 + \lambda_1 X_1 + \lambda_2 \E[X_2|\upto{1}{l}] - \bX_1 \psi_1 (\aopt{1}-a_1).
\end{eqnarray*}
It is clear from this form that if $\E[X_2|\upto{1}{l}]$ depends on $a_1$, then the modelled blip function $\gammaopt{1}(a_1 ,\upto{1}{l};\psi_1)$ does not capture the entire effect of $A_1$ (that is, there is mediation of the effect of $A_1$ on $Y$ through $X_2$).  This problem is essentially identical to that described in the sections above describing direct forward simulation.
\index{Regret|)}
\medskip

\noindent \textbf{Simulation:}  We use a Monte Carlo study of 5000 replicate data sets, and study the statistical properties of the estimators of the blip parameters.  The treatment-free outcome models for the two intervals were specified as
\[
\mopt{02}  :  \quad \beta_{02} + \beta_{12} a_1 + \beta_{22} x_1 + \beta_{32} x_2 \qquad \qquad \mopt{01}   : \quad \beta_{01} + \beta_{11} x_1
\]
but are clearly mis-specified compared to the true data generating model.  The results are contained in Table \ref{tab:Ex4c}; even for very large sample sizes, there is a degree of bias especially for the first stage parameters, and the rate at which the optimal regime is correctly identified is only around 81 \% even when the sample size is 10000.

\begin{table}[ht]
\caption{Example \ref{ex:Ex5desc}: Proportion accuracy for identifying the true optimal regime, $\sqrt{n}  \times \textrm{Bias}$ and $n \times \textrm{Variance}$ of G-estimators computed using 5000 replicate data sets and five sample sizes under a correctly specified exposure model.}\label{tab:Ex4c}
  \centering
 \begin{tabular}{|lr|r|r|r|r|r|r|}
    \hline
&&\multicolumn{5}{c|}{Sample size $n$}\\[3pt]
&&500&1000&2000&5000&10000\\
\hline
\multicolumn{2}{|r|}{Proportion accuracy}&0.509&0.583&0.660&0.753&0.812\\[3pt]
\hline
$\sqrt{n} \times \textrm{Bias}$ & $\psi_{11}$ &9.100&20.777&34.953&62.921&93.243
\\
&$\psi_{21}$& 46.127&63.177&90.231&141.532&198.480\\
&$\psi_{31}$& -7.966&-11.876&-18.373&-31.670&-45.372\\
\cline{2-7}
&$\psi_{12}$& 20.387&30.093&43.716&69.885&98.876\\
&$\psi_{22}$& 18.420&16.781&15.766&18.415&23.893\\
&$\psi_{32}$& 15.620&17.827&25.048&37.093&48.865
\\
&$\psi_{42}$&-1.771&-8.008&-14.205&-25.720&-37.977\\
&$\psi_{52}$&18.095&-15.109&-15.713&-18.447&-22.570\\
\hline
$n \times \textrm{Var}$ & $\psi_{11}$ &976.83&586.959&436.83&363.558&328.795
\\
&$\psi_{21}$& 3506.467&2005.572&1712.518&1601.715&1635.553
\\
&$\psi_{31}$&1501.585&1080.018&972.281&952.681&947.784
\\
\cline{2-7}
&$\psi_{12}$& 787.487&516.049&458.349&443.737&415.136
\\
&$\psi_{22}$& 5756.547&2814.746&2099.608&1843.867&1692.712
\\
&$\psi_{32}$& 5272.212&2890.574&2486.51&2279.34&2100.907
\\
&$\psi_{42}$& 879.778&490.78&397.897&340.086&337.187
\\
&$\psi_{52}$& 1539.997&863.299&691.944&607.725&584.84
\\
    \hline
  \end{tabular}
  \vspace{0.1 in}
\end{table}

\subsection{\texorpdfstring{$K$}{K} intervals}

The direct simulation approach for $K$ intervals proceeds in principle as above, and considers a chain rule generation for confounder and treatment patterns across time, specifically the conditional models
\[
X_k \; | \; \upto{k-1}{A}, \upto{k-1}{X} \qquad \textrm{and} \qquad A_k \; | \; \upto{k-1}{A}, \upto{k}{X}
\]
to generate complete profiles $(\upto{K}{A},\upto{K}{X})$.  We would then consider the complete conditional model
$Y | \upto{K}{A}, \upto{k}{X}$ through its expectation; for example, we might specify
\[
\E[Y|\upto{K}{A} = \upto{K}{a}, \upto{K}{X} = \upto{K}{x}] = \left( \lambda_0 + \sum_{k=1}^K \lambda_{1k} x_k + \sum_{k=1}^K \sum_{l=1}^K \lambda_{2kl} x_k x_l \right) + \left(\sum_{k=1}^K \lambda_{3k} a_k + \sum_{k=1}^K \lambda_{4k} a_k x_k \right)
\]
for parameters $\lambda$.  This induces a collection of parameters $\psi_k,k=1,\ldots,K$ that parameterize the true (data generating) optimal blip-to-zero models $\gammaopt{k}(.,.;
\psi_k), k=1,\ldots,K$.

It is evident from the previous discussion, however, that such an approach will require certain structural assumptions (for example, equality constraints on the $\lambda$ parameters), due to the need for the optimal blip-to-zero function\index{Blip function} $\gammaopt{k}(.,.;
\psi_k)$ to be functionally independent of $\aopt{k+1},\ldots,\aopt{K}$ for each $k$.  Therefore, simulation using the regret-based strategy described in section \ref{sec:Regretbased} is in general more straightforward.  The drawback of the regret-based method is again that the treatment free outcome models are not generally straightforward to specify correctly, and thus the robustness to mis-specification of the G-estimation method is important.

\begin{myexample} \label{ex:Ex6Desc}
For $K=5$, with binary treatment, suppose that we specify optimal blip-to-zero models
\[
\gammaopt{k}(a,\upto{k}{h};\psi_k) = a (\psi_{1k} + \psi_{2k} x_k + \psi_{3k} a_{k-1} + \psi_{4k} x_k a_{k-1}) \qquad k=1,2,\ldots,5,
\]
with $a_0 \equiv 0$, and with the model
\[
\E[Y(\uptoopt{K}{a})|\upto{K}{X}=\upto{K}{x}] = \lambda_0 + \lambda_1 x_1
\]
characterizing the outcome under optimal treatment.  For the treatment model, suppose that
\[
\Pr[A_k = 1 | \upto{k}{H} = \upto{k}{h}] = \textrm{expit}\{\alpha_{0k} + \alpha_{1k} x_k + \alpha_{2k} a_{k-1} \}.
\]
We fit correctly specified treatment and optimal blip-to-zero models in each interval, and specify a treatment-free outcome model of the form
\[
\mopt{0k}(\upto{k}{h};\bbeta_k) = \beta_{0k} + \beta_{1k} a_{k-1} + \beta_{2k} x_{k-1} + \beta_{3k} x_k \qquad k=1,\ldots,K
\]
where $x_{0} \equiv 0$.

\medskip

\noindent \textbf{Simulation:}  A Monte Carlo study of 5000 replicate data sets with the following parameter settings was carried out:
\begin{eqnarray*}
\lambda & : \qquad  & (0.2,-0.1)\\
\sigma_X & : \qquad & 1\\
\alpha_k & : \qquad & (0.5,0.5,-0.2), \quad k=1,\ldots,K\\
\psi_k & : \qquad & (0.5,1.0), k=1, \;(-1.0,0.5,2.0,-0.2), \quad k=2,\ldots,K\\
\sigma_Y & : \qquad & 1.
\end{eqnarray*}
Five different sample sizes $n=500,1000,2000,5000,10000$ were used.

\medskip
The results, presented as scaled bias and variance, are presented in Tables \ref{tab:Ex6m} and \ref{tab:Ex6v}.  Both scaled bias and variance increase with $n$ over this range, and also are higher for the later intervals.  However, the proportion of times that the correct optimal regime is identified increases as $n$ increases: even at $n=500$, the correct regime is identified with around 80
\% accuracy.

\medskip
\noindent Code for this example is available on request from the author.

\begin{table}[ht]
\caption{Example \ref{ex:Ex6Desc}: Proportion accuracy for identifying the true optimal regime, $\sqrt{n}  \times \textrm{Bias}$ of G-estimators computed using 5000 replicate data sets and five sample sizes under a correctly specified exposure model.}\label{tab:Ex6m}
  \centering
 \begin{tabular}{|c|r|r|r|r|r|r|}
    \hline
&\multicolumn{5}{c|}{Sample size $n$}\\[3pt]
&500&1000&2000&5000&10000\\
\hline
\multicolumn{1}{|r|}{Acc.}&0.809&0.866&0.901&0.934&0.953\\
\hline
$\psi_{11}$&  8.153& 13.370& 20.185& 33.956& 48.813\\
$\psi_{21}$& 20.578& 30.130& 43.314& 69.476& 98.837\\
\hline
$\psi_{12}$&-22.759&-31.766&-44.802&-70.223&-99.414\\
$\psi_{22}$& 10.532& 15.203& 21.678& 34.839& 49.477\\
$\psi_{32}$& 41.484& 60.520& 86.908&139.080&197.707\\
$\psi_{42}$& -4.284& -6.072& -8.651&-13.973&-19.740\\
\hline
$\psi_{13}$&-22.386&-31.383&-44.520&-69.995&-99.127\\
$\psi_{23}$& 10.525& 15.290& 21.945& 35.030& 49.432\\
$\psi_{33}$& 41.792& 60.819& 87.237&139.099&197.641\\
$\psi_{43}$& -4.145& -6.215& -9.024&-14.149&-19.681\\
\hline
$\psi_{14}$&-22.087&-31.026&-44.239&-69.842&-99.026\\
$\psi_{24}$& 10.713& 15.348& 21.872& 34.885& 49.407\\
$\psi_{34}$& 42.501& 60.933& 87.342&139.267&197.686\\
$\psi_{44}$& -4.400& -6.227& -8.874&-13.982&-19.603\\
\hline
$\psi_{15}$&-21.271&-30.381&-43.798&-69.599&-98.892\\
$\psi_{25}$& 10.587& 15.056& 21.757& 34.908& 49.621\\
$\psi_{35}$& 42.620& 61.011& 87.519&139.476&198.015\\
$\psi_{45}$& -4.262& -5.899& -8.622&-14.072&-20.004\\
    \hline
  \end{tabular}
  \vspace{0.1 in}
\end{table}

\begin{table}[ht]
\caption{Example \ref{ex:Ex6Desc}: $n \times \textrm{Variance}$ of G-estimators computed using 5000 replicate data sets and five sample sizes under a correctly specified exposure model.}\label{tab:Ex6v}
  \centering
 \begin{tabular}{|c|r|r|r|r|r|r|}
    \hline
&\multicolumn{5}{c|}{Sample size $n$}\\[3pt]
&500&1000&2000&5000&10000\\
\hline
$\psi_{11}$& 54.075& 55.098& 59.153& 61.965& 60.877\\
$\psi_{21}$&  6.577&  6.435&  6.108&  5.927&  6.109\\
\hline
$\psi_{12}$& 73.892& 71.446& 75.225& 79.640& 77.259\\
$\psi_{22}$& 26.637& 23.743& 24.134& 24.086& 22.752\\
$\psi_{32}$& 50.828& 47.595& 47.413& 47.652& 46.594\\
$\psi_{42}$& 36.792& 33.574& 34.097& 33.541& 31.739\\
\hline
$\psi_{13}$& 90.276& 85.959& 92.482& 93.869& 93.829\\
$\psi_{23}$& 57.482& 56.426& 53.609& 57.252& 54.245\\
$\psi_{33}$&105.198&102.444&100.715& 97.058&101.648\\
$\psi_{43}$& 81.428& 76.446& 74.047& 76.855& 75.675\\
\hline
$\psi_{14}$& 73.986& 74.595& 78.389& 77.120& 74.722\\
$\psi_{24}$& 60.323& 56.259& 55.843& 53.625& 53.756\\
$\psi_{34}$&113.560&111.944&112.915&112.820&109.270\\
$\psi_{44}$& 86.811& 82.362& 81.666& 79.354& 78.343\\
\hline
$\psi_{15}$& 62.347& 60.915& 59.384& 58.254& 59.832\\
$\psi_{25}$& 64.879& 64.397& 62.127& 58.802& 61.031\\
$\psi_{35}$&129.385&126.639&126.217&121.716&129.979\\
$\psi_{45}$& 93.475& 91.699& 91.732& 87.903& 89.749\\
    \hline
  \end{tabular}
  \vspace{0.1 in}
\end{table}

\end{myexample}

\section{Summary and Discussion}

In this paper we have seen formulation of the inference problem for an optimal dynamic treatment regime based on G-estimation.  The development of a hyperopic treatment strategy which prioritizes long-term benefits over short-term gains is based on a backwards induction\index{Backward induction} approach: in the $K$ intervals of treatment setting, the optimal treatment at interval $k \leq K$ is inferred under the assumption that optimal treatment is allocated for intervals $l > k$.

The G-estimation procedure is based on a semiparametric model, applied sequentially across intervals, that at interval $k$ specifies the conditional expectation, given history up to interval $k$, of a pseudo-outcome that is derived from the final outcome by subtracting the estimated effects of exposure at intervals $k+1,k+2,\ldots,K$.   The conditional expectation is partitioned into two components, the first capturing the treatment-free contribution, the second the additional contribution due to sub-optimal treatment, at that interval.  The G-estimating equation utilizes also the model for treatment at interval $k$ given history up to that interval in the construction of a doubly robust\index{Doubly robust} estimator that is consistent for the parameters in the optimal blip-to-zero treatment outcome model if at least one of the treatment-free outcome model and the treatment allocation model is correctly specified.   Given the consistency in estimation, the optimal treatment pattern across intervals may also be consistently reconstructed; the optimal dynamic treatment regime is determined via a partitioning of the space of subject histories which constructs a series of allocation rules across the intervals.\index{Blip function}

G-estimation can be viewed as an extended form of least-squares estimation where the estimating equation is amended using the treatment allocation model.  Under standard regularity conditions, the resulting estimators possess properties -- consistency and asymptotic normality -- determined by the usual semiparametric frequentist theory.  However, assumptions made about the optimal blip-to-zero treatment model
\[
\gammaopt{k}(\upto{k-1}{a},a,\upto{k}{H};\psi_k) = \E[Y(\upto{k-1}{a},a,\afteropt{k+1}{a})-Y(\upto{k-1}{a},0,\afteropt{k+2}{a})|\upto{k}{H};\psi_k] \qquad k=1,\ldots,K
\]
are key: all inferences, and treatment decisions, are dependent on this structural assumption.  It is therefore important to verify that this model component is correctly specified.

In this paper, the comparison of outcomes under hyperopic regimes with those derived under myopic strategies has not been carried out explicitly.  It is easy to verify that the hyperopic strategy produces the optimal dynamic regimen: see for example \cite{Rich2013,Rich2015}.
\index{G-estimation|)}

\vskip .2in
\noindent{\bf Acknowledgements}

\noindent This research was supported by a Discovery Grant from the Natural Sciences and Engineering Research Council of Canada (NSERC).

\bibliographystyle{plain}
\bibliography{DTRbook-References}

\end{document}